\documentclass{scrartcl}

\usepackage{graphicx}
\graphicspath{{./images/}}
\usepackage{amsmath,amsthm,amssymb}
\usepackage[numbers]{natbib}
\usepackage{xr-hyper}
\usepackage[hypertexnames=false,colorlinks=true,linkcolor=blue,citecolor=blue,anchorcolor=blue]{hyperref}
\usepackage[capitalise]{cleveref}
\usepackage{color}
\usepackage{acro}
\usepackage{mathtools}
\usepackage{booktabs}
\usepackage{pdfpages}

\definecolor{mylightblue}{RGB}{73, 223, 252}

\newcommand{\quark}{\setbox0\hbox{$x$}\hbox to\wd0{\hss$\cdot$\hss}}

\newcommand{\email}[1]{\protect\href{mailto:#1}{#1}}

\newcommand\funding[1]{\protect\\ \hspace*{1.8em}{\bfseries Funding:} #1}
\newcommand{\mylozenge}{$\mathbin{\blacklozenge}$}
\newcommand{\R}{\mathbb{R}}

\DeclareMathOperator{\tr}{tr}

\newtheorem{theorem}{Theorem}

\newtheorem{proposition}{Proposition}

\theoremstyle{definition}

\newtheorem{assumption}{Assumption}
\crefname{assumption}{Assumption}{Assumptions}

\theoremstyle{remark}
\newtheorem{remark}{\textbf{Remark}}  
\crefname{remark}{Remark}{Remarks}

\theoremstyle{definition}

\DeclareAcronym{statfem}{short=StatFEM, long=statistical finite element method}
\DeclareAcronym{fem}{short=FEM, long=Finite Element Method}
\DeclareAcronym{pde}{short=PDE, long=partial differential equation}
\DeclareAcronym{ode}{short=ODE, long=ordinary differential equation}
\DeclareAcronym{rkhs}{short=RKHS, long=reproducing kernel Hilbert space}
\DeclareAcronym{gp}{short=GP, long=Gaussian process,short-plural=s,long-plural=es}
\DeclareAcronym{bvp}{short=BVP, long=boundary value problem}
\DeclareAcronym{femesh}{short=FE mesh, long=finite element mesh,short-plural=es,long-plural=es}
\DeclareAcronym{dof}{short=DOF, long=degree of freedom,long-plural-form=degrees of freedom}
\DeclareAcronym{fe}{short=FE, long=Finite Element}
\DeclareAcronym{pn}{short=PN, long=Probabilistic Numerics}
\DeclareAcronym{pnm}{short=PNM, long=Probabilistic Numerical Method}
\DeclareAcronym{ivp}{short=IVP, long=initial value problem}
\DeclareAcronym{gmrf}{short=GMRF, long=Gaussian Markov Random Field}

\author{Yanni Papandreou\thanks{Imperial College London (\email{john.papandreou18@imperial.ac.uk})} \and Jon Cockayne \thanks{The Alan Turing Institute (\email{jcockayne@turing.ac.uk})} \and Mark Girolami \thanks{The Alan Turing Institute and University of Cambridge (\email{mag92@cam.ac.uk})} \and Andrew B. Duncan\thanks{Imperial College London and The Alan Turing Institute (\email{a.duncan@imperial.ac.uk})}}
\title{Theoretical Guarantees for the Statistical Finite Element Method\thanks{\funding{This work was supported by Wave 1 of The UKRI Strategic Priorities Fund under the EPSRC Grant EP/T001569/1 and EPSRC Grant EP/W006022/1, particularly the “Ecosystems of Digital Twins” theme within those grants \& The Alan Turing Institute. YP was supported by a Roth Scholarship funded by the Department of Mathematics, Imperial College London.}}}

\ifpdf
\hypersetup{
  pdftitle={Theoretical Guarantees for the Statistical Finite Element Method},
  pdfauthor={Y. Papandreou, J. Cockayne, M. Girolami, and A. B. Duncan}
}
\fi

\begin{document}

\maketitle

\begin{abstract}
    The statistical finite element method (StatFEM) is an emerging probabilistic method that allows observations of a physical system to be synthesised with the numerical solution of a PDE intended to describe it in a coherent statistical framework, to compensate for model error.
    This work presents a new theoretical analysis of the statistical finite element method demonstrating that it has similar convergence properties to the finite element method on which it is based.
    Our results constitute a bound on the Wasserstein-2 distance between the ideal prior and posterior and the StatFEM approximation thereof, and show that this distance converges at the same mesh-dependent rate as finite element solutions converge to the true solution.
    Several numerical examples are presented to demonstrate our theory, including an example which test the robustness of StatFEM when extended to \emph{nonlinear} quantities of interest.
\end{abstract}

\section{Introduction}
\label{sec:intro}

Mathematical models based on \acp{pde} play a central role in science and engineering. Since PDEs do not generally yield closed form solutions, approximate solutions are typically obtained via numerical methods. The \ac{fem} \cite{strang1973analysis} is one such scheme which approximates the solution by solving a discretised problem over a finite mesh of the PDE domain. \ac{fem}-based methods have been widely studied, and as such, rigorous estimates of the error introduced through discretisation have been obtained for various families of PDEs. In particular, these results characterise the convergence of the error to zero as a function of the size of the mesh.

Mismatch between the model and the real system can occur for a variety of reasons, for example, if the geometry of the domain is not fully resolved, if underlying processes are not completely understood, or if approximations must be introduced to reduce the model to a closed system of PDEs. Stochasticity is often introduced in models to reflect the uncertainty arising from this lack of knowledge, or simply to capture intrinsically noisy processes. In the context of PDEs, stochasticity may appear in various forms, in the initial and/or boundary conditions, in the forcing terms, in the coefficients of the PDE or even in the geometry of the domain on which the PDE is defined. The resulting governing equations, known as Stochastic Partial Differential Equations (SPDEs), can be reformulated probabilistically as a probability distribution over the space of admissible solutions. The development and analysis of numerical methods for SPDEs is a very active area of research, see \cite{zhang2017numerical} for a review. 

The increasing prevalence of instrumentation within  modern engineering systems has driven the need for tools able to combine sensor data with predictions from PDE models in a manner which can mitigate any model misspecification. This technology forms an integral part of the modern notion of a digital twin, as a means of providing effective continuous monitoring of complex engineering assets. 

The recently proposed \ac{statfem} framework \cite{girolami2019statistical,Girolami2021} seeks to address this need through a  statistical formulation of \ac{fem}, which permits numerical approximations of a PDE to be combined with sensor observations in a principled manner through the Bayesian formalism.  Due to its general applicability and ability to provide quantification of uncertainty without expensive Monte Carlo simulations, it has established itself as an effective means of assimilating data into the large-scale PDEs which arise in civil, marine and aerospace engineering. As \ac{statfem} is built on \ac{fem} it is itself subject to numerical discretisation error which will propagate through to any predictions made by the model. Obtaining rigorous bounds on the bias introduced by discretisation error on \ac{statfem} predictions and their associated uncertainty is a critical gap in the literature that this work will seek to address.

\subsection{Related Work}
\label{sec:related}

This paper is focused on providing detailed error analysis for the \ac{statfem} method first introduced in \cite{girolami2019statistical,Girolami2021}. For clarity of ideas, we limit our attention to linear elliptic problems. Adapting these results to time-dependent and/or non-linear PDEs is left to future research.

There has been extensive research on deriving numerical methods for solving various classes of SPDEs.  For detailed reviews of such methods the reader is referred to \cite{lord2014introduction,MATTHIES1997283,sudret2000stochastic,STEFANOU20091031,xiu2010numerical,aldosary2018structural}. Some of these methods, for example the Stochastic Galerkin Finite Element method (SG-FEM) and its many variants, also finds their roots in the classical finite element method, but differ fundamentally from StatFEM. More specifically, SG-FEM methods reformulate the SPDE as a deterministic PDE on a higher dimensional domain, where additional variables are introduced to characterise the underlying noise process. On the other hand, in the StatFEM approach, the approximation is itself intrinsically stochastic.   While SG-FEM methods are applicable to a wider class of SPDEs, their scalability is limited by the curse of dimensionality resulting from the augmented noise dimensions. Additionally, while StatFEM formulation permits fast assimilation of measurement data, this is not straightforward for SG-FEM and related methods.

In the area of spatial statistics, Gaussian random fields are constructed as the solution of appropriately chosen SPDEs \cite{lindgren2011explicit}. The solution is approximated through a finite element discretisation of the weak form of the underlying SPDE. The resulting discretised solution is subsequently employed as a prior within a Bayesian inference problem for a spatial or spatio-temporal model \cite{sigrist2015stochastic}. For non-Gaussian likelihoods, fast approximate Bayesian Inference for the resulting latent Gaussian model is possible using an Integrated Nested Laplace Approximation (INLA) \cite{martino2014integrated}. The construction of the \ac{statfem} model is strongly connected to this statistical viewpoint of SPDEs.  Indeed, in both cases, the spatial priors are \acp{gmrf}. In the case of \ac{statfem} the mean and covariance of the \ac{gmrf} reflect the behaviour of the underlying deterministic PDE, while in the SPDE case the covariance is tailored to ensure a specific degree of spatial decorrelation of the underlying latent Gaussian random field.

Conceptually, \ac{statfem} also bears some similarity to Probabilistic Numerics (PN) Methods \cite{hennig2015probabilistic}, both being statistical reformulations of numerical methods.  Indeed, PN methods have also been developed for approximating solutions of \acp{ode} \cite{tronarp2019probabilistic} as well as PDEs \cite{conrad2017statistical,chkrebtii2016bayesian,cockayne2016probabilistic}.   The two methodologies share common foundations, particularly Bayesian Probabilistic Numerical (BPN) Methods \cite{cockayne2019bayesian} which are formulated as conditioning a prior distribution of potential solutions on partial evaluations of the underlying equations.  \ac{statfem} also builds on the Bayesian formalism, employing conditioning to incorporate new information, specifically noisy sensor observations, starting from a prior which characterises the distribution of approximate PDE solutions.  However, these methods fundamentally differ in what they seek to achieve:   \acp{pnm} methods seek to recast numerical discretisation error as uncertainty,  whereas the uncertainty arising in \ac{statfem} primarily seeks to account for model error.

\subsection{Contributions}

The novel contributions of the paper are as follows:
\begin{itemize}
	\item We study the bias introduced by discretisation error in the \ac{statfem} model, and derive rigorous upper bounds for the bias introduced in the prior and posterior in terms of the properties of the underlying finite element mesh.
	\item We conduct a thorough simulation study confirming our theoretical results on a number of illustrative \ac{fem} test problems.
	\item We empirically explore the limits of our theory by studying the bias induced by \ac{statfem} in a scenario where the assumptions underpinning our theoretical results are no longer valid. 
\end{itemize}

\subsection{Structure of the Paper}

The remainder of the paper is structured as follows. In \cref{sec:notation} we introduce notation required for the rest of the paper. In \cref{sec:statfem} we reintroduce \ac{statfem}, deriving it from first principles. \Cref{sec:theory} contains our main theoretical contributions, error bounds for the \ac{statfem} prior and the posterior when this prior is conditioned on sensor observations. We conduct simulations verifying our theoretical results in \cref{sec:experiments}, and providing some concluding remarks in \cref{sec:conclusion}. The supplementary material contains the proofs required for the paper in \cref*{supp-sec:supplement}.

\subsection{Notation} \label{sec:notation}

In this section we introduce required notation for the remainder of the paper.  Given a normed space $\mathcal{X}$ denote the associated norm by $\|\cdot\|_{\mathcal{X}}$. Given two normed spaces $\mathcal{X},\mathcal{Y}$ we  denote by $L(\mathcal{X},\mathcal{Y})$ the set of all bounded linear operators from $\mathcal{X}$ to $\mathcal{Y}$. We denote the dual of $\mathcal{X}$ by $\mathcal{X}^{\prime}=L(\mathcal{X},\R)$.  When $\mathcal{X}$ is a Hilbert space we  denote the associated inner product by $\langle \cdot, \cdot \rangle_{\mathcal{X}}$. We will sometimes drop the subscript for the norm and inner product when there is no risk of confusion over which spaces we are referring to.

Let $\mathcal{X},\mathcal{Y}$ be Hilbert spaces. Given an operator $A\in L(\mathcal{X},\mathcal{Y})$ we denote the adjoint of $A$ by $A^{\dagger}\in L(\mathcal{Y},\mathcal{X})$. Given $A\in L(\mathcal{X},\mathcal{X})$, we say that $A$ is symmetric if $A = A^\dagger$ and positive definite if it is symmetric and $\langle x, Ax\rangle > 0$ for all nonzero  $x\in \mathcal{X}$. Recall that the trace of $A$ can be defined as $\tr(A)=\sum_{i=1}^{\infty}\langle Ae_i,e_i \rangle$ where $\left\{ e_i \right\}_{i=1}^{\infty}$ is any orthonormal basis of $\mathcal{X}$.  We will say that $A\in L(\mathcal{X},\mathcal{X})$ is trace-class if $\tr(A) < \infty$.  For an operator $A:\mathcal{X}\rightarrow \mathcal{Y}$ we will define the operator norm to be $\|A\|_{\mathcal{X}\rightarrow\mathcal{Y}}=\sup_{x\in \mathcal{X}} \frac{\|Ax\|_{\mathcal{Y}}}{\|x\|_{\mathcal{X}}}$, the trace norm to be $\|A\|_{\tr}=\tr([A^{\dagger}A]^{1/2})$, while the Hilbert-Schmidt norm is defined to be $\|A\|_{\text{HS}}=\tr(A^{\dagger}A)^{1/2}$. Recall the well known inequality $\|A\|_{\mathcal{X}\rightarrow\mathcal{Y}}\leq\|A\|_{\text{HS}}\leq\|A\|_{\tr}$.

Let $D\subset\R^{d}$ be an open domain with boundary $\partial D$ and let $\overline{D} = D \cup \partial D$ be its closure. We denote by $C(D), C(\overline{D})$ the space of continuous functions on $D$ and $\overline{D}$, respectively. Denote by $L^p(D)$ the spaces of $p$-integrable functions on $D$ and $W^{k,p}(D)$ the Sobolev space with $p$-integrable derivatives up to and including order $k$.  We use the notation $H^k(D) = W^{k,2}(D)$ and $H_{0}^k(D)$ to be the elements of $H^k(D)$ which are zero on $\partial D$. For clarity we will often use the notation $\|\cdot\|_{k}$, $|\cdot|_{k}$ to denote the norm and semi-norm associated to $H^{k}(D)$.

\subsection{Relevant Probability-Theoretic Concepts}
\label{sec:prob_concepts}
In this section we briefly introduce and discuss the relevant probability-theoretic concepts we will require. 
To this end, for a space $\mathcal{X}$ let $\mathcal{B}(\mathcal{X})$ denote the Borel sigma-algebra of $\mathcal{X}$, and for a measurable space $(\mathcal{X}, \mathcal{B}(\mathcal{X}))$ let $\mathcal{P}(\mathcal{X})$ denote the set of all Borel probability measures on $\mathcal{X}$.

Given a positive definite kernel function $k:D\times D\rightarrow \mathbb{R}$ on the domain $D$, and a function $m:D\rightarrow \mathbb{R}$ we say $f$ is a Gaussian process with mean function
$m$ and covariance function $k$ if for every finite collection of points $\lbrace x_1,\ldots, x_N\rbrace \subset D$ the random vector $(f(x_1), . . . , f(x_N))$ is a multivariate Gaussian random variable with mean vector $(m(x_1), . . . , m(x_N))$ and covariance matrix $\left\{ k(x_n, x_m) \right\}_{n,m=1,\ldots, N}$.  The mean function and kernel/covariance function completely determine the Gaussian process. We write
$f \sim \mathcal{GP}(m, k)$ to denote the Gaussian process with mean function $m$ and covariance function $k$.  We refer the reader to \cite{Rasmussen06gaussianprocesses} for further details. 

Gaussian processes that take values in a separable Hilbert space $\mathcal{X}$ can be associated canonically with a Gaussian measure on the space $\mathcal{X}$.  Gaussian measures on infinite-dimensional Hilbert spaces can be viewed as a generalisation of their finite-dimensional counterpart on $\mathbb{R}^d$, and are  defined by a mean element $m\in \mathcal{X}$ and covariance operator $K$, which is a positive-definite trace class operator. See \cite{da2014stochastic} and \cite{bogachev1998gaussian} for more information.  Given a Gaussian process $f \sim \mathcal{GP}(m, k)$, this can be associated with the Gaussian measure $\mathcal{N}(m, K)$ with mean $m$ and covariance operator $K$, associated with the kernel $k$, defined by
$$
Kg(\cdot) = \int k(\cdot, y)g(y)\,dy, \quad g \in \mathcal{X}.
$$
Given any $m \in \mathcal{X}$ and $K \in \mathcal{L}(\mathcal{X},\mathcal{X})$ which is positive-definite and trace class, there exists a Gaussian measure $\mathcal{N}(m,K)$ with mean $m$ and covariance operator $K$ \cite[Proposition 2.18]{da2014stochastic}. 

For the purposes of studying the error in \ac{statfem} we will need to introduce a metric on the space of distributions. The probability metric we will work with in this paper is the Wasserstein distance; this will be justified in \cref{sec:theory}. Given two probability measures $\mu,\nu$ on a normed space $\mathcal{X}$, having the first two moments finite, the Wasserstein-2 distance between the two measures, $W_2(\mu,\nu)$, is defined as:
\begin{equation}
	\label{def:wasser}
	W_2^{2}(\mu,\nu)=\inf_{\pi\in\Gamma(\mu,\nu)}\int_{\mathcal{X}\times\mathcal{X}}\|x-y\|_{\mathcal{X}}^{2}\mathrm{d}\pi(x,y)
\end{equation}
where $\Gamma(\mu,\nu)$ is the set of couplings of $\mu$ and $\nu$, i.e.
\begin{align*}
	\Gamma(\mu,\nu)=\{ \pi \in \mathcal{P}(\mathcal{X} \times \mathcal{X}) \mid \pi(E\times\mathcal{X})=\mu(E)
	\text{ and } \pi(\mathcal{X}\times F)=\nu(F) \text{ for all } E,F\in \mathcal{B}(\mathcal{X}) \}.
\end{align*}
Throughout we will abbreviate $W_2(\mu, \nu) = W(\mu, \nu)$.

It will be also useful to note that the Wasserstein distance between the measures $\mu,\nu$ above can be linked to the Wasserstein distance between the centered measures.
Let $m$ denote the expectation of $\mu$; then the centered measure $\mu^*$ is defined by
\begin{equation*}
	\mu^*(E) = \mu(\{x - m \mid x \in E\})
\end{equation*}
for all $E \in \mathcal{B}(\mathcal{X})$. 
Denoting the means of $\mu,\nu$ by $m_{1},m_{2}$ respectively we have from \cite{cuesta1996lower}:
\begin{equation}\label{eq:wasser_centred}
	W^{2}(\mu,\nu)=\|m_{1}-m_{2}\|_{\mathcal{X}}^{2}+W^{2}(\mu^{*},\nu^{*}).
\end{equation}
We also recall that there is an explicit expression for the Wasserstein distance between two Gaussians. Let $\mu_{i}=\mathcal{N}(m_{i},\Sigma_{i}), i=1,2$ be two Gaussian measures on a separable Hilbert space $\mathcal{X}$. Then from \cite{masarotto2019procrustes}, we have
\begin{equation}\label{eq:wasser_for_gaussians}
	W^{2}(\mu_{1},\mu_{2}) = \|m_1-m_2\|_{\mathcal{X}}^{2}+\tr(\Sigma_{1})+\tr(\Sigma_{2})-2\tr\sqrt{\Sigma_{1}^{1/2}\Sigma_{2}\Sigma_{1}^{1/2}} .
\end{equation}  

\section{The Statistical Finite Element Method} 
\label{sec:statfem}
In this section we will present the \ac{statfem} model. We begin by reviewing \ac{fem} in \cref{sec:fem}. Then, in \cref{sec:statfem_prior} we will derive the prior distribution over the solution to the \ac{pde} that is used in \ac{statfem}.
Finally, in \cref{sec:statfem_posterior} we will condition this prior distribution on sensor observations to obtain the \ac{statfem} posterior distribution.

\subsection{The Finite Element Method}
\label{sec:fem}

In this section we will briefly review the necessary concepts from \ac{fem} needed for \ac{statfem}, focusing on elliptic \acp{pde}. We begin by establishing notation. Consider the following linear elliptic \ac{bvp}:

\begin{equation}\label{eq:elliptic_bvp}
	\begin{aligned}
	\mathcal{L}u(x):=	-\nabla\cdot\left( \kappa(x)\nabla u(x) \right) &= f(x) &\quad x &\in D \\
		u(x) &= 0 &\quad x&\in\partial D,
	\end{aligned}
\end{equation}
for some positive, spatially-dependent conductivity  $\kappa:D\rightarrow\R_{> 0}$ and forcing $f\in L^{2}(D)$. We make the following assumptions on the the \emph{data} $(D, \kappa, f)$:

\begin{assumption}[domain regularity]\label{ass:domain_reg}
	The dimension of the domain is restricted to $d\in\left\{ 1,2 \right\}$. The domain $D$ is convex polygonal, i.e.\ it is convex and its boundary $\partial D$ is a polygon. Note that this implies that the domain is bounded. \mylozenge
\end{assumption}

\Cref{ass:domain_reg} is assumed for simplicity; we note that it should be possible to relax some of these assumptions though this is not investigated in this work.

\begin{assumption}[parameter regularity]
	\label{ass:param_reg}
	The parameter $\kappa$ is a continuous function on $\overline{D}$, i.e.\ $\kappa\in C(\overline{D})$. Further, there exists a constant $\alpha>0$ such that:
	\begin{equation}\label{eq:param_ellipticity}
		\kappa(x)\geq\alpha > 0, \quad \text{for almost all } x\in D
	\end{equation}
	Inequality \cref{eq:param_ellipticity} is often referred to as \textit{(uniform) ellipticity} of the parameter. \mylozenge
\end{assumption}

We consider the \textit{weak solution} of the \ac{bvp} \cref{eq:elliptic_bvp}, namely we seek $u\in H_{0}^{1}(D)$ which satisfies the variational formulation of \cref{eq:elliptic_bvp}, i.e.

\begin{equation}
	\label{eq:variation_form}
	\int_{D}\nabla v(x)\cdot\left( \kappa(x)\nabla u(x) \right) \mathrm{d}x=\int_{D}f(x)v(x)\mathrm{d}x, \quad \forall v\in H_{0}^{1}(D)
\end{equation}

Under \cref{ass:domain_reg,ass:param_reg}, standard regularity results \cite[section 6.2.2]{evans10} show that, for each $f\in L^{2}(D)$, there exists a unique \textit{weak solution} $u\in H_{0}^{1}(D)$ satisfying \cref{eq:variation_form}.
Further, \cref{ass:domain_reg,ass:param_reg} guarantee stronger regularity of the weak solution. 
From \cite[section 3.7]{larsson2008partial}  it holds that $u\in H^{2}(D)$ and that there exists a constant $\tilde{C}>0$ such that
$\|u\|_{2}\leq \tilde{C}\|f\|_{L^{2}(D)}$. Note that this extra regularity can be generalised to elliptic problems in higher dimensions \cite[Theorem 3.2.1.2]{grisvard2011elliptic}

We can reformulate the variational formulation \cref{eq:variation_form} as 
\begin{equation}
	\label{eq:shortened_var_form}
	a(u,v)=\langle f,v \rangle \quad \forall v\in H_{0}^{1}(D)
\end{equation}
where the bilinear form $a:H_{0}^{1}(D)\times H_{0}^{1}(D)\rightarrow\R$ is defined by 
\begin{equation}\label{eq:def_of_bilinear_form}
	a(u,v):=\int_{D}\nabla v(x)\cdot \left(\kappa(x)\nabla u(x)\right) \mathrm{d}x
\end{equation} 
\cref{ass:param_reg} on $\kappa$ ensures that the bilinear form $a$ is symmetric, bounded and coercive \cite{evans10}.

The \emph{Galerkin approximation} of the variational problem \eqref{eq:shortened_var_form} seeks an approximate solution $u_h \in V_h$ within a finite dimensional subspace $V_h$ of $H_0^1(D)$ which satisfies
$$
    a(u_h, v) = \langle f, v\rangle, \quad \forall v\in V_h.
$$ The construction of the finite element method arises from a specific Galerkin approximation, where we introduce a \textit{triangulation} of the domain, $\mathcal{T}_{h}=\left\{ \omega_e \right\}_{e=1}^{n_e}$, where the $\omega_e$ are pairwise disjoint triangular elements with $D = \cup_{e=1}^{n_e} \omega_e$, and $h$ is the diameter of the largest element. This triangulation is usually referred to as the \ac{femesh} and we will refer to $h$ as the \textit{mesh width}.  The discrete function space $V_h$ is then chosen to be piecewise polynomial functions defined on this mesh.  The approximate weak solution $u_h$ can be expressed as a linear combination of basis functions $\lbrace \phi_i \rbrace_{i=1}^{n_u}$.  In this work, we assume that $V_h$ is the space of piecewise continuous functions, and $\lbrace \phi_i \rbrace_{i=1}^{n_u}$ chosen to be a \emph{nodal basis} which satisfies $\phi_i(x_j) = \delta_{ij}$, where $x_1,\ldots, x_{n_u}$ are the internal vertices of the mesh.

By exploiting the linearity of both the bilinear form and the inner product in \cref{eq:shortened_var_form}
the \ac{fem} discretisation of the variational formulation of our \ac{bvp} yields the following discrete system of linear equations:
\begin{equation}
	\label{eq:fem_system}
	A\mathbf{u}=\mathbf{f}
\end{equation}
where $A\in\R^{n_u\times n_u}$ is the \textit{stiffness matrix} whose $ij$\textit{-th} element is given by:
\begin{equation}
	\label{eq:stiffness_mat}
	A_{ij}=a(\phi_j,\phi_i)=\int_{D}\nabla\phi_{i}\cdot(\kappa(x)\nabla\phi_j)\mathrm{d}x
\end{equation}
The vector $\mathbf{u}\in\R^{n_u}$ is the vector of \acp{dof} and $\mathbf{f}\in\R^{n_u}$ is the \textit{load vector} whose $i$\textit{-th} entry is given by:
\begin{equation}
	\label{eq:load_vect}
	f_{i}=\langle f,\phi_{i} \rangle=\int_{D}f\phi_{i}\mathrm{d}x
\end{equation}
The \textit{stiffness matrix}, $A$, is a real symmetric matrix which is invertible by \cref{ass:param_reg} \cite[see Theorem 2.16]{lord2014introduction}. Thus, the system of linear equations \cref{eq:fem_system} can be solved for $\mathbf{u}$, yielding the finite element approximation  $u_h=\boldsymbol{\Phi}^{T}\mathbf{u}$ where $\boldsymbol{\Phi}:=\left( \phi_1,\cdots,\phi_{n_u} \right)^{T}$ is the vector of the FE basis functions. The \textit{mass matrix}, $M\in\R^{n_u\times n_u}$  which is defined by
\begin{equation}
	\label{eq:mass_mat}
	M_{ij}=\langle \phi_i,\phi_j \rangle=\int_{D}\phi_{i}\phi_{j}\mathrm{d}x,
\end{equation}
will also play an important role in the error analysis of \cref{sec:theory}.

It will be important for our work to have access to a bound on the error between the true solution and the \ac{fem} approximation in terms of $h$. To this end consider a family of triangulations $\{ \mathcal{T}_{h_i} \}$, such that $h_i \to 0$ and $h_i > 0$ for all $i$.
We will require a further technical assumption on the family of the \acp{femesh}:

\begin{assumption}
	\label{ass:refinement_reg}
	The meshes under consideration remain regular in the sense that as we decrease $h_{i}$ to $0$ the angles of all triangles are bounded below independently of $h_{i}$. \mylozenge
\end{assumption}

Under the assumptions above we have the following classical $H^{1}$ error bound \cite[see chapter 5]{larsson2008partial}:
\begin{theorem}
	\label{thm:fem_error_H1}
	Let $u$ be the unique \textit{weak solution} to \ac{bvp} \cref{eq:elliptic_bvp} and let $u_h$ by the \ac{fem} approximation. 
	Assume that $f\in L^{2}(D)$ and that the \cref{ass:domain_reg,ass:param_reg,ass:refinement_reg} hold. 
	We then have the following $H^{1}$ error bound, for sufficiently small $h$:
	\begin{equation}
		\label{eq:fem_error_bound_H1}
		|u-u_h|_{1}\leq Ch\|u\|_{2}
	\end{equation}
	where $C>0$ is a positive constant independent of $h$. \mylozenge
\end{theorem}

From this $H^{1}$ error bound a standard duality argument yields the following $L^{2}$ error bound \cite[see chapter 5]{larsson2008partial}, which will play an important role in the error analysis of \cref{sec:theory}.

\begin{theorem}
	\label{thm:fem_error}
	Let $u$ be the unique \textit{weak solution} to \ac{bvp} \cref{eq:elliptic_bvp} and let $u_h$ by the \ac{fem} approximation. 
	Assume that $f\in L^{2}(D)$ and that the \cref{ass:domain_reg,ass:param_reg,ass:refinement_reg} hold. 
	We then have the following error bound, for sufficiently small $h$:
	\begin{equation}
		\label{eq:fem_error_bound}
		\|u-u_h\|_{L^{2}(D)}\leq Ch^{2}\|f\|_{L^{2}(D)}
	\end{equation}
	where $C>0$ is a positive constant independent of $h$. \mylozenge
\end{theorem}

\subsection{The \ac{statfem} Prior} \label{sec:statfem_prior}

In this section we present the \ac{statfem} approximation for the solution of a stochastically-forced elliptic PDE. The resulting probability distribution will subsequently be employed as a prior to be combined with sensor data through Bayesian updates.  To this end, we consider the \ac{bvp} \cref{eq:elliptic_bvp} where $f$ is chosen to be a \ac{gp},
$$
f \sim \mathcal{GP}(\bar{f}, k_{f}),
$$
where $\bar{f} \in C(D)$ and where we assume the covariance $k_f$ is regular enough so that realisations of $f$ lie almost surely in $L^{2}(D)$. It is sufficient to take $k_f(x,y)=\psi(x-y)$ where $\psi:D\rightarrow\R$ is continuous \cite[see Corollary 1]{scheuerer2010regularity}. This is satisfied for example by the RBF kernel or any kernel in the Mat\'ern family. Another alternative sufficient condition would be to choose the kernel so that the associated covariance operator is an appropriate negative power of the Laplacian. This is outlined in detail in \cite[see section 6]{stuart2010inverse}.

In our analysis we shall assume that the conductivity coefficient $\kappa$ is deterministic, unlike \cite{Girolami2021} where this is assumed to be stochastic. Since \cref{eq:elliptic_bvp} is well-posed, the solution can be expressed as $u=\mathcal{L}^{-1}f$, where $\mathcal{L}^{-1}$ is the inverse of the differential operator $\mathcal{L}$. The pushforward of the probability measure for $f$ through $\mathcal{L}^{-1}$ induces a probability distribution of solutions supported on $H_0^1(D)$.  Since the inverse operator is linear, the resulting distribution for $u$ will also be a Gaussian measure supported on $L^2(D)$, with mean $\mathcal{L}^{-1}\overline{f}$ and covariance operator $\mathcal{L}^{-1} K \mathcal{L}^{-\dagger}$, where $K$ is the covariance operator of $f$ on $L^2(D)$. Indeed, $u$ can be expressed as a Gaussian process \cite[Proposition 3.1]{owhadi2015bayesian}
\begin{equation}
    \label{eq:true_distr}
    u \sim \mathcal{GP}(\mathcal{L}^{-1}\bar{f}, k_{\mathcal{L}})
\end{equation}
where
$$
    k_{\mathcal{L}}(x,y) = \mathcal{L}_x^{-1}\mathcal{L}_y^{-1} k_f(x,y)
$$

Here $\mathcal{L}_x^{-1}$ denotes the action of the operator $\mathcal{L}^{-1}$ on the first argument of the kernel and $\mathcal{L}^{-1}_y$ its action on the second argument.

We intend to use \cref{eq:true_distr} as a prior, which will later be combined with sensor observations through sequential Bayesian updates. An important observation is that the mean and covariance of the GP defined in \cref{eq:true_distr} cannot be computed exactly for general problems, which precludes the calculation of marginal distributions or realisations of the GP. To overcome this issue, another GP, which employs a finite element approximation of the mean and covariance functions will be introduced. 

As discussed in \cref{sec:fem} \ac{fem} approximates the solution to our elliptic \ac{bvp} as a linear combination of the FE basis functions. In particular, the \ac{fem} approximation to the solution of our elliptic \ac{bvp} can be expressed as follows:
\begin{equation}
	u_{h}=\boldsymbol{\Phi}^{T}A^{-1}\mathbf{f}=:\hat{L}^{-1}_{h}f
\end{equation}
where the operator $\hat{L}_{h}^{-1}:L^{2}(D)\rightarrow L^{2}(D)$ is defined as $\hat{L}_{h}^{-1}:=\boldsymbol{\Phi}^{T}A^{-1}\mathcal{F}_{h}$. Here the operator $\mathcal{F}_{h}:L^{2}(D)\rightarrow\R^{n_{u}}$ takes a function $f$ to the load vector $\mathbf{f}$, defined in \cref{eq:load_vect}, and $A$ is the stiffness matrix defined in \cref{eq:stiffness_mat}. The \ac{statfem} approximation to \eqref{eq:true_distr}  is obtained by formally replacing  the solution operator $\mathcal{L}^{-1}$ with $\hat{L}_{h}^{-1}$ in the mean and covariance terms.  Due to the linearity of the approximation solution map, the resulting process is still a Gaussian process, with transformed mean and covariance. 

\begin{proposition}
	\label{prop:statfem_prior}
		The \ac{statfem} prior from \cite{girolami2019statistical,Girolami2021} is a \ac{gp} of the form
		\begin{equation}
			\label{eq:statfem_prior_gp}
			u_{h}\sim\mathcal{GP}(\bar{u}_{h},k_h)
		\end{equation}
		where the mean and covariance function are defined by:
		\begin{align}
			\bar{u}_{h}(x) &:= \boldsymbol{\Phi}(x)^{T}A^{-1}\bar{\mathbf{f}} \\
			k_{h}(x,y) &:= \boldsymbol{\Phi}(x)^{T}A^{-1}K_{\mathcal{F}_{h}}A^{-1}\boldsymbol{\Phi}(y)
		\end{align}
		where $K_{\mathcal{F}_{h}}:=\mathcal{F}_{h}K\mathcal{F}_{h}^{\dagger}\in\R^{n_u\times n_u}$ and $\bar{\mathbf{f}}:=\mathcal{F}_{h}\bar{f}\in\R^{n_u}$. \mylozenge
	\end{proposition}

The question of how much error is introduced by replacing $u$ with $u_h$ can be formulating as quantifying the distance between the probability measures associated with each Gaussian process. To this end, we define $\nu_\star$ and $\nu_h$ to be the Gaussian probability measures supported on $L^2(D)$ which are associated with the GPs \eqref{eq:true_distr} and \eqref{eq:statfem_prior_gp}, respectively.  We can write $\nu_\star = \mathcal{N}(m_\star, \Sigma_\star)$ and $\nu_h = \mathcal{N}(m_h, \Sigma_h)$ where 
\begin{align*}
	m_{\star} &= \mathcal{L}^{-1}\bar{f} \\
	\Sigma_{\star} &= \mathcal{L}^{-1}K\mathcal{L}^{-\dagger},
\end{align*}
and
\begin{align*}
	m_{h} &= \hat{L}_{h}^{-1}\bar{f} \\
	\Sigma_{h} &= \hat{L}_{h}^{-1}K\hat{L}_{h}^{-\dagger}
\end{align*}

\begin{remark}
	\label{rem:justification_of_mutual_singularity}
	Let $V_{\star}, V_{h}$ be the supports of $\nu_{\star}, \nu_{h}$ respectively. Since the image of the approximate solution operator $\hat{L}_{h}^{-1}$ is the span of the FE basis functions $\left\{ \phi_{i} \right\}_{i=1}^{n_{u}}$ it follows that, for all $h>0$, $V_h$ is a finite dimensional subspace of $H_{0}^{1}(D)\subseteq L^{2}(D)$. Further, under \cref{ass:domain_reg,ass:param_reg} we have that $V_{\star}$ is an infinite dimensional subspace of $H_{0}^{1}(D)\subseteq L^{2}(D)$. In the case that $V_h$ and $V_{\star}$ are disjoint we clearly have that $\nu_{\star}$ and $\nu_{h}$ are mutually singular. Consider now the case that $V_h$ and $V_{\star}$ are not disjoint and $V_h\not\subseteq V_{\star}$. This means that $V_h$ and $V_{\star}$ have a non-empty intersection, and that there exists a subset $\widetilde{V}\subseteq V_{\star}$ not in the intersection. Further, $\nu_{h}$ gives 0 measure to $\widetilde{V}$, while $\widetilde{V}$ has positive measure under $\nu_{\star}$. This implies that $\nu_{\star}$ and $\nu_{h}$ are mutually singular. Finally, in the case that $V_{h}\subset V_{\star}$ we can write $V_{\star}$ as the orthogonal direct sum of $V_{h}$ and its complement. Since $\nu_{h}$ gives full measure to $V_h$ it must give measure 0 to the orthogonal complement of $V_h$, which has positive measure under $\nu_{\star}$. Hence, we again have that $\nu_{\star}$ and $\nu_{h}$ are mutually singular. In all cases, the mutual singularity follows from the Feldman-Hajek theorem \cite[Theorem 2.25]{da2014stochastic}.
\end{remark}

\subsection{Posterior Updates of the  \ac{statfem} approximation} \label{sec:statfem_posterior}

We now consider the problem of inferring the solution $u$ of the BVP \eqref{eq:elliptic_bvp} for an unknown forcing term $f$ based on noisy observations of the form 
\begin{equation}
	\label{eq:sensor_likelihood}
	\mathbf{v} = Su + \boldsymbol{\xi}, \quad \boldsymbol{\xi} \sim \mathcal{N}(0, \epsilon^2 I)
\end{equation}
where $S$ is a linear observation map.   Following the Bayesian formalism \cite{stuart2010inverse}, the solution of this inverse problem is characterised by the posterior distribution $\widetilde{\nu}_\star(du)$ of $u$ given the observations $\mathbf{v}$,
\begin{equation}
\label{eq:posterior_exact}
   \widetilde{\nu}_\star(du) \propto  e^{-\Phi(\mathbf{v}; u)}  {\nu}_\star(du)
\end{equation}
where $\Phi(\mathbf{v}; u) = \frac{1}{2\epsilon^2 }\lVert Su - \mathbf{v}\lVert^2,
$
is the negative log-likelihood, and $\epsilon$ is the observational noise standard deviation.   In practice, we would work with an alternative posterior distribution, where the prior distribution $\nu_{\star}$ is replaced with a prior obtained from a finite element approximation. To understand the consequences of making this approximation, we must quantify the distance (in an appropriate sense) between the posterior probability measures \eqref{eq:posterior_exact} and 
\begin{equation}
\label{eq:posterior_approximate}
    \widetilde{\nu}_h(du) \propto  e^{-\Phi(\mathbf{v}; u)}  {\nu}_h(du).
\end{equation}
We shall make the following assumptions, on the observation map $S$. 

\begin{assumption} 
\label{ass:observation}The observation map $S:C(D)\rightarrow \mathbb{R}^s$ is a bounded linear operator with full range.
\end{assumption}

\begin{remark}
	\label{rem:domain_of_obs_op}
	Note that the assumptions on the covariance structure of the forcing $f$ guarantee that solutions to our noisy elliptic \ac{bvp} lie in $H^{2}(D)$ almost surely. Thus, the support of the true prior lies in $H^{2}(D)$. By the Sobolev Embedding Theorem for dimensions 1,2, one has that $H^{2}(D)\subset C(D)$ so that $S$ is well-defined for functions drawn from the true prior. Draws from the \ac{statfem} prior are trivially in $C(D)$ and so $S$ is well-defined for functions drawn from this prior as well. Note that the condition on the range of $S$ will always be true for some $s > 0$.
\end{remark}

A property which will be very important for the subsequent error analysis of Section \ref{sec:theory} is that the posterior  distributions $\widetilde{\nu}_\star$ and $\widetilde{\nu}_h$ remain Gaussian measures. This can be viewed as the infinite dimensional version of the analogous result for multivariate Gaussians based on ``completing the square''. 

\begin{proposition} \label{prop:statfem_posteriors} Suppose that the observation map satisfies Assumption \ref{ass:observation} and let $\mathbf{v} \in \mathbb{R}^s$.  Then $\widetilde{\nu}_\star$ and $\widetilde{\nu}_h$ are Gaussian measures with mean and covariance given by 

\begin{align}
	\label{eq:post_means}
	m_{u|\mathbf{v}}^{(a)}&:=m_{a} + \Sigma_{a}S^{\dagger}(\epsilon^{2}I+S\Sigma_{a}S^{\dagger})^{-1}(\mathbf{v}-Sm_{a}) \\
	\label{eq:post_covs}
	\Sigma_{u|\mathbf{v}}^{(a)}&:=\Sigma_{a} - \Sigma_{a}S^{\dagger}(\epsilon^{2}I+S\Sigma_{a}S^{\dagger})^{-1}S\Sigma_{a},
\end{align}
respectively, for $a \in \lbrace \star, h\rbrace$.\textbf{}
\mylozenge
\end{proposition}

\section{Discretisation Error Analysis} 
\label{sec:theory}

Intuitively we expect that $\nu_{h}\rightarrow\nu_{\star}$ as $h\rightarrow 0$ in an appropriate sense because of the convergence properties of the finite element discretisation. Similarly, we expect that the approximate posterior, $\tilde{\nu}_{h}$, will converge to the true posterior, $\tilde{\nu}_{\star}$, in an appropriate sense, as $h\rightarrow 0$. 

To establish these properties we must select an appropriate metric for the distance between the priors and posteriors. A challenge with selecting this metric is that, as explained in \cref{rem:justification_of_mutual_singularity}, the measures $\{\nu_{h},\nu_{\star}\}$ and $\{\tilde{\nu}_{h},\tilde{\nu}_{\star}\}$ are generally mutually singular. This complicates the typical stability analysis that arises when studying Bayesian inverse problems, based on Kullback-Leibler divergence or  Hellinger distance \cite{cotter2009bayesian}. In this paper we shall adopt an alternative approach, employing the Wasserstein-2 metric to quantify the influence of discretisation error on the posterior distribution. This is advantageous in this context for several reasons: firstly it is not dependent on the use of a common dominating measure. Secondly, unlike other divergences, it is robust under vanishing noise: in this limit, the Wasserstein distance between the two probability measures will converge to the deterministic $L^2$ discretisation error between the respective two means. Finally, while computing Wasserstein distance is intractable for general measures, there exists a closed form for the Wasserstein-2 distance between Gaussian measures, valid in both finite and infinite dimensions. See \cite{masarotto2019procrustes} for a statement of this and see \cite{dowson1982frechet,olkin1982distance} for proofs in the finite dimensional case, and \cite{cuesta1996lower} for a proof in the infinite dimensional case.

The results in this section will focus on bounding the Wasserstein distance between the priors $\nu_\star$ and $\nu_h$, and the posteriors $\widetilde{\nu}_\star$ and $\widetilde{\nu}_h$ as a function of $h$. The strategy employed to obtain these bounds involves exploiting an important connection between Wasserstein distance between Gaussian measures on Hilbert spaces and the Procrustes Metric between the respective covariance operators, as detailed in \cite{masarotto2019procrustes}. We now present these results in the following two theorems.

\begin{theorem}[Prior Error Analysis]
    \label{thm:statfem_prior_error}
    Assuming that $f\in L^{2}(D)$, and that the domain, $D$, and parameter, $\kappa$, satisfy \cref{ass:domain_reg,ass:param_reg,ass:refinement_reg} there exists a constant $\gamma>0$, independent of $h$, such that, for sufficiently small $h$ we have:
    \begin{equation}
        \label{eq:statfem_prior_error}
        W(\nu_{\star},\nu_{h})\leq\gamma h^{2}
    \end{equation} \mylozenge
\end{theorem}

The next theorem describes how this error in the prior propagates forward to the error between arbitrary linear functionals of the true and \ac{statfem} posterior, again as a function of the mesh width $h$.

\begin{theorem}[Posterior Error Analysis]
    \label{thm:statfem_posterior_error}
    Let $\ell:L^{2}(D)\rightarrow\R$ be a bounded linear functional. Assuming that $f\in L^{2}(D)$, and that the domain, $D$, and parameter, $\kappa$, and the observation map $S$, satisfy \cref{ass:domain_reg,ass:param_reg,ass:refinement_reg,ass:observation} there exists a constant $\gamma^{\prime}(\ell)>0$, independent of $h$, such that, for sufficiently small $h$ we have:
    \begin{equation}
        \label{eq:statfem_posterior_error}
        W(\ell_{\#}\tilde{\nu}_{\star},\ell_{\#}\tilde{\nu}_{h})\leq\gamma^{\prime}(\ell)h^{2} + \mathcal{O}(h^{4}) \text{ as } h\rightarrow 0
    \end{equation}
    \mylozenge
\end{theorem}

\begin{remark}
    \cref{thm:statfem_posterior_error} provides a bound on $W(\ell_{\#}\tilde{\nu}_\star, \ell_{\#}\tilde{\nu}_h)$, for arbitrary $\ell$, rather than on $W(\tilde{\nu}_\star, \tilde{\nu}_h)$.
    Naturally a bound on the latter quantity is stronger and more of interest, but deriving such a bound is technically challenging.
    Nevertheless the empirical results in \cref{sec:experiments} suggest that the stronger result may hold, and as a result this will be the focus of future research.
    \mylozenge
\end{remark}

\cref{thm:statfem_prior_error} shows that we have control over the Wasserstein distance between the full priors while \cref{thm:statfem_posterior_error} shows that, in the posterior case, we have control over linear quantities of interest. Nonlinear quantities of interest are also interesting and we will explore whether they are controlled empirically in the numerical experiments in \cref{sec:experiments}.

\section{Experiments} 
\label{sec:experiments}

In this section we will present some numerical experiments demonstrating the theoretical results from \cref{sec:theory}. In \cref{sec:1D_poisson} we will consider a one-dimensional example in which we have access to an explicit formula for the Green's function, while in \cref{sec:2D_poisson} we will consider a more challenging two-dimensional problem in which the solution operator is not available. In both scenarios we will aim to demonstrate the \ac{statfem} prior and posterior bounds.
Note that while \cref{thm:statfem_posterior_error} shows only that $W(\ell_{\#}\tilde{\nu}_\star, \ell_{\#}\tilde{\nu}_h)$ is $\mathcal{O}(h^2)$ our experiments suggest that for the problems examined here it also holds that $W(\tilde{\nu}_\star, \tilde{\nu}_h)$ is $\mathcal{O}(h^2)$. 
Finally, in \cref{sec:1D_maximum} we explore whether our convergence results extend beyond the linear setting by calculating the Wasserstein distance between the prior and posterior distributions of a \emph{nonlinear} quantity of interest, given by $\max_{x\in D}u(x)$. The code for these numerical experiments is available online\footnote{\url{https://github.com/YanniPapandreou/statFEM}}. Additional implementation details for the experiments are provided in \cref*{supp-sec:impl_details}.

\subsection{One-dimensional Poisson Equation}\label{sec:1D_poisson}
In this section we will present a numerical experiment involving the following one-dimensional Poisson equation:
\begin{equation}\label{eq:1D_poisson_bvp}
    \begin{aligned}
        -&\frac{d^{2}u}{dx^{2}}(x) = f(x) \quad x \in [0,1] \\
        &u(0)=u(1)=0
    \end{aligned}
\end{equation}
where the forcing $f$ is distributed as,
\begin{equation}\label{eq:1D_forcing_distr}
    f \sim \mathcal{GP}(\bar{f},k_{f}) \\
\end{equation}
with mean and covariance function given by:
\begin{align}
    \label{eq:1D_forcing_distr_mean}
    \bar{f}(x) &= 1 \\
    \label{eq:1D_forcing_distr_cov}
    k_{f}(x,y) &= \sigma_{f}^{2}\exp\left( -\frac{|x-y|^{2}}{2l_{f}^{2}} \right).
\end{align}
We took $\sigma_f=0.1$ and $l_f=0.4$. Note that the \ac{bvp} \cref{eq:1D_poisson_bvp} is an instance of the general problem \cref{eq:elliptic_bvp} with $d=1$, $D=[0,1]$, and $\kappa(x)\equiv 1$. The Green's function, $G(x,y)$, for the \ac{bvp} \cref{eq:1D_poisson_bvp} is easily computed to be:
\begin{equation}\label{eq:1D_greens_func}
    G(x,y) = x(1-y)\Theta(y-x) + (1-x)y\Theta(x-y) \quad\forall x,y\in[0,1]
\end{equation}
where $\Theta(x)$ is the Heaviside Step function. As discussed in \cref{sec:statfem_prior} the solution $u$ to problem \cref{eq:1D_poisson_bvp} will be the Gaussian process defined by \cref{eq:true_distr}. Since we readily have the Green's function available for our 1D problem we can explicitly write down this distribution as:
\begin{equation}\label{eq:1D_sol_distr}
    u\sim\mathcal{N}(\mu_{\star},k_{\star})
\end{equation}
where the mean and covariance functions are given by:
\begin{align}
    \label{eq:1D_sol_mean}
    \mu_{\star}(x)&=\int_{0}^{1}G(x,w)\bar{f}(w)\mathrm{d}w=\frac{1}{2}x(1-x), \\
    k_{\star}(x,y)&=\int_{0}^{1}\int_{0}^{1}G(x,w)k_{f}(w,t)G(t,y)\mathrm{d}t\mathrm{d}w
\end{align}
Note that we have explicitly computed the integral for the mean $\mu_{\star}$, but not for the covariance $k_{\star}$. The double integral for the covariance is approximated numerically using quadrature to accurately compute this on a reference grid.

\subsubsection{Prior results for 1-D example}\label{sec:1D_poisson_prior}
Since we can straightforwardly evaluate the true prior numerically, we can now demonstrate the error bound for the \ac{statfem} prior given in \cref{thm:statfem_prior_error}. In order to do this the following remark, taken from \cite{convergenceRateSite}, will be useful.
\begin{remark}\label{rem:conv_rate_true_res_known}
    When we have access to the true distributions we can estimate the rate of convergence of the \ac{statfem} priors and posteriors as follows. Let $\eta_{\star},\eta_{h}$ denote the true and \ac{statfem} prior (or posterior) respectively. If we assume that we have a rate of convergence $p$, i.e.,
    \begin{equation}\label{eq:conv_rate_assump}
        W(\eta_{\star},\eta_{h})\leq\mathcal{O}(h^{p}) \text{ as } h\rightarrow 0
    \end{equation}
    then we have that there is a constant $C>0$, independent of $h$, such that,
    \begin{equation}\label{eq:conv_rate_assump_impl}
        W(\eta_{\star},\eta_{h})=Ch^{p} + \mathcal{O}(h^{p+1})
    \end{equation}
    from which it follows that,
    \begin{equation}
        \log W(\eta_{\star},\eta_{h}) = p\log h + \log|C| + \mathcal{O}(h)
    \end{equation}
    From this we can see that we can estimate the rate of convergence, $p$, when we know the true distributions by computing the Wasserstein distance for a range of sufficiently small $h$-values, plotting a log-log plot of these values, and estimating the slope of the line of best fit.
    \mylozenge
\end{remark}

As outlined in \cref{rem:conv_rate_true_res_known} above we compute, for a range of sufficiently small $h$-values, an approximation of the Wasserstein distance between the two priors and then plot these results on a log-log scale. We take a reference grid of $N=51$ equally spaced points in $[0,1]$ to compare the covariances, and we take $30$ $h$-values in $[0.02,0.25]$. With these choices we obtain the plot shown in \cref{fig:1D_prior_results} below. From this plot we can see that the results indeed lie along a straight line with a slope estimated to be $p=2.00$ (to $2$ decimal places) using \texttt{scipy}'s built in linear regression function.

\begin{figure}[!ht]
    \centering
    \includegraphics[width=12cm]{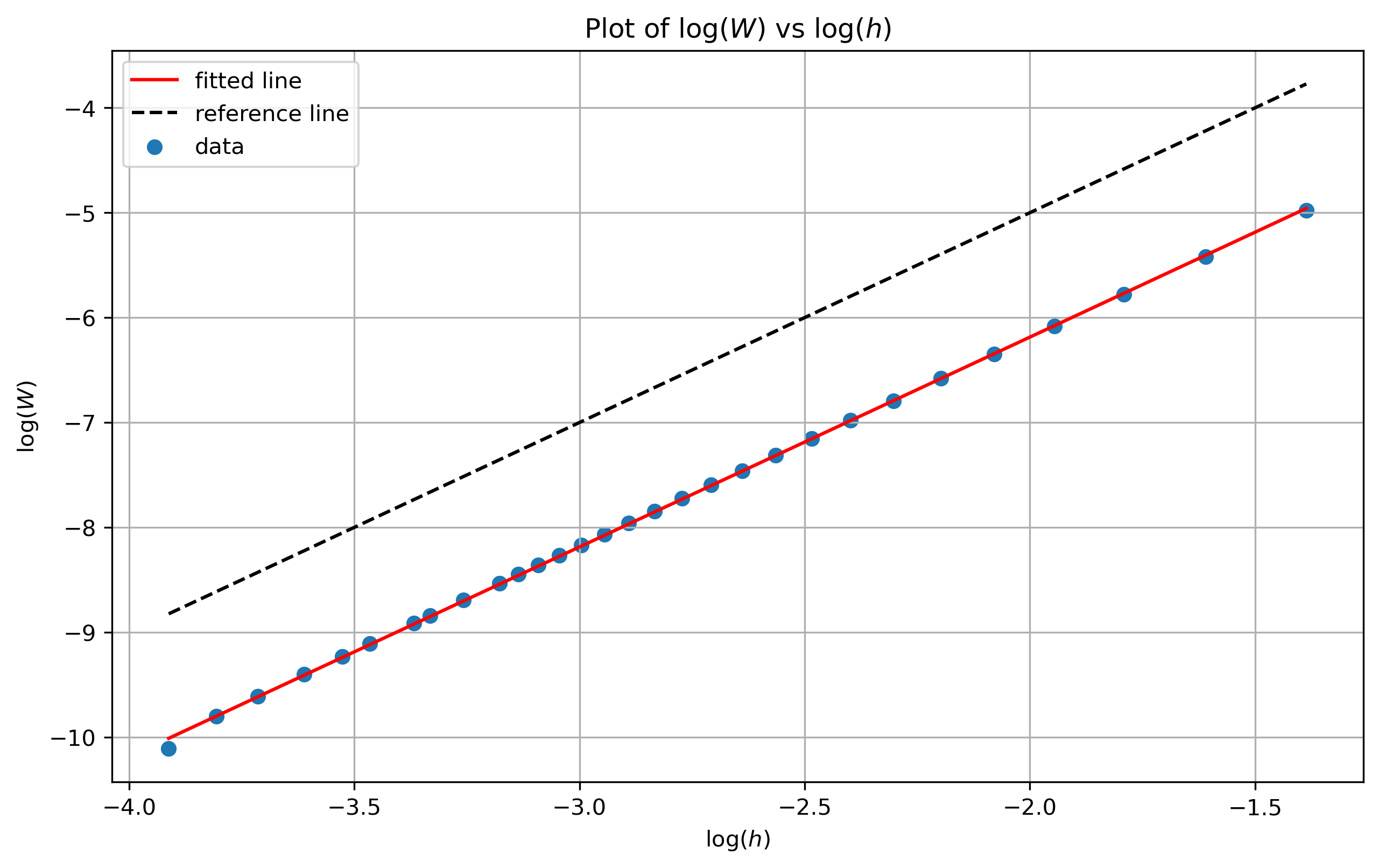}
    \caption{Log-log plot of Wasserstein distance against $h$ for the priors in a 1-D example. The reference line has a slope of 2. The slope of the fitted line is estimated to be $p=2.00$, in agreement visually with the reference line.}
    \label{fig:1D_prior_results}
\end{figure}

\subsubsection{Posterior results for 1-D example}\label{sec:1D_poisson_posterior}
For this example we will consider a specific case of the general Bayesian inverse problem outlined in \cref{sec:statfem_posterior} which will correspond to taking the observation operator $S$ to be a pointwise evaluation operator. More specifically, we suppose that we have sensor data which correspond to noisy observations of the value of the solution $u$ at some points $\left\{ y_{i} \right\}_{i=1}^{s}\subset D$. We thus take $S$ to be the operator which maps a function $g$ to the vector $Sg:=(g(y_1),\dots,g(y_s))^{T}\in\R^{s}$. For this choice of $S$ the update equations for the posterior means and covariances given in \cref{eq:post_means,eq:post_covs} respectively can be simplified to a form which is more amenable to computation. These forms, for $i$ being either the symbol $\star$ or $h$, are:
\begin{align}
    \label{eq:post_means_eval_form}
    m_{u|\mathbf{v}}^{(i)}(x) &= m_{i}(x)-\mathbf{k}^{(i)}(x)^{T}B_{\epsilon,i}^{-1}(\mathbf{m}^{(i)}-\mathbf{v}) \\
    \label{eq:post_covs_eval_form}
    k_{u|\mathbf{v}}^{(i)}(x,y) &= k_{i}(x,y) - \mathbf{k}^{(i)}(x)^{T}B_{\epsilon,i}^{-1}\mathbf{k}^{(i)}(y) \\
    \label{eq:B_mat_form}
    B_{\epsilon,i} &:= \epsilon^{2}I+S\Sigma_{i}S^{\dagger}=\epsilon^{2}I+C_{Y,i} \\
    \label{eq:cov_vect_form}
    \mathbf{k}^{(i)}(x) &:= (k_{i}(x,y_{1}),\dots,k_{i}(x,y_{s}))^{T}
\end{align}
where $k_{u|\mathbf{v}}^{(i)}$ is the covariance function associated with $\Sigma_{u|\mathbf{v}}^{(i)}$ and where the matrix $C_{Y,i}:=S\Sigma_{i}S^{\dagger}\in\R^{s\times s}$ has $pq$\textit{-th} entry $k_{i}(y_{p},y_{q})$.
The specific sensor observations used in the experiments were obtained by simulating a trajectory from $\nu_\star$ evaluated at the locations $y_1, \dots, y_s$.

Since we have the true posterior in an explicit form we can now demonstrate that the error bound for the full \ac{statfem}
posterior agrees with the error bound for linear quantities of interest from the posterior given in \cref{thm:statfem_posterior_error}. We demonstrate this error bound for 4 different levels of sensor noise, $\epsilon=0.00005,0.0001,0.01,0.1$. These levels are chosen so that the lower values are comparable to the prior variances. We take $s=10$ sensor readings equally spaced in the interval $[0.01,0.99]\subset D$. For each level of sensor noise we again follow the argument outlined in \cref{rem:conv_rate_true_res_known} and compute, for a range of sufficiently small $h$-values, an approximation of the Wasserstein distance between the two posteriors, and plot these results on a log-log scale. We take a reference grid of $N=41$ equally spaced points in $[0,1]$ to compare the covariances, and we take 28 $h$-values in $[0.025,0.25]$. With these choices we obtain the plot shown in \cref{fig:1D_posterior_results} below. We also present the estimates for the slopes and intercepts of the lines of best fit in \cref{tab:1D_posterior_results}.

\begin{figure}[!ht]
    \centering
    \includegraphics[width=12cm]{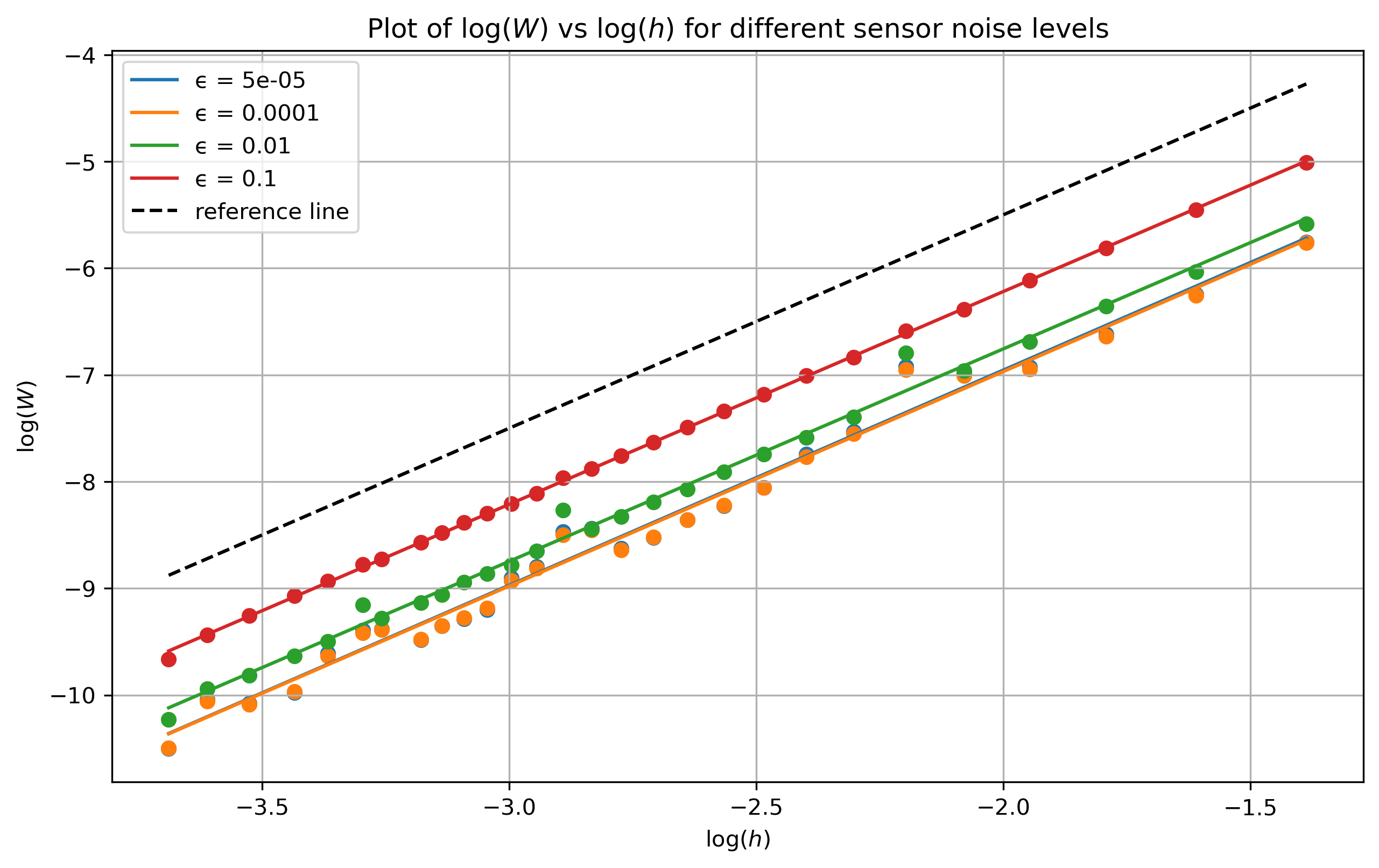}
    \caption{Log-log plot of Wasserstein distance against $h$ for the posteriors in a 1-D example for varying levels of sensor noise. The reference line has a slope of 2. Notice how the best fit lines are all parallel to the reference line.}
    \label{fig:1D_posterior_results}
\end{figure}

\begin{table}[tbhp]
    {\footnotesize
    \caption{Estimate of slopes and intercepts for the different levels of sensor noise $\epsilon$. All of the slopes are around 2.}\label{tab:1D_posterior_results}
    \begin{center}
        \begin{tabular}{@{}ccc@{}}
            \toprule
            $\epsilon$ & slope  & intercept \\ \midrule
            0.00005       & 2.0162 & -2.9227   \\
            0.00010       & 2.0102 & -2.9486   \\
            0.01000       & 1.9912 & -2.7739   \\
            0.10000       & 1.9940 & -2.2309   \\ \bottomrule
        \end{tabular}
    \end{center}
    }
\end{table}

We can see from \cref{fig:1D_posterior_results} and \cref{tab:1D_posterior_results} that we obtain lines of best-fit which are all parallel to each other, each with a slope of around $p=2$. The intercepts are slightly different for each level of sensor noise, as is to be expected since the constant in \cref{thm:statfem_posterior_error} depends upon $\epsilon$. These results show that the error bound for the full posterior distributions is also of order $p=2$, suggesting that it is possible to extend \cref{thm:statfem_posterior_error} to the Wasserstein distance between the full posteriors. Proving this is left to future research.

\subsection{Two-dimensional Poisson Equation}\label{sec:2D_poisson}

In this section we will present a numerical example involving the following two-dimensional Poisson equation:
\begin{equation}\label{eq:2D_poisson_bvp}
	\begin{aligned}
		-\Delta u(x) &= f(x) &\quad x &\in D:=[0,1]^{2} \\
		u(x) &= 0 &\quad x&\in\partial D
	\end{aligned}
\end{equation}
where the forcing $f$ is distributed as,
\begin{equation}\label{eq:2D_forcing_distr}
    f \sim \mathcal{GP}(\bar{f},k_{f}) \\
\end{equation}
with mean and covariance given by
\begin{align}
    \label{eq:2D_forcing_distr_mean}
    \bar{f}(x) &= 1 \\
    \label{eq:2D_forcing_distr_cov}
    k_{f}(x,y) &= \sigma_{f}^{2}\exp\left( -\frac{\|x-y\|^{2}}{2l_{f}^{2}} \right).
\end{align}
Again we took $\sigma_{f}=0.1$ and $l_{f}=0.4$. Note that the \ac{bvp} \cref{eq:2D_poisson_bvp} is an instance of the general problem \cref{eq:elliptic_bvp} with $d=2$, $D=[0,1]^{2}$, and $\kappa(x)\equiv 1$.  For this problem we use the approach outlined in \cref{rem:conv_rate_true_res_unknown} below to demonstrate the error bounds. We present the error estimates for the prior and posterior separately in the following sections.

\subsubsection{Prior results for 2-D example}\label{sec:2D_poisson_prior}

Following the notation of \cref{thm:statfem_prior_error} we denote the \ac{statfem} prior by $\nu_{h}$. Since we do not have access to the true distributions here we will need an alternative approach to that of \cref{rem:conv_rate_true_res_known}. This approach is outlined in the remark, taken from \cite{convergenceRateSite}, below:
\begin{remark}\label{rem:conv_rate_true_res_unknown}
    When we do not have access to the true distributions we can proceed as follows. Again assuming that we have a rate of convergence, $p$, as given by \cref{eq:conv_rate_assump}, we can obtain an estimate as follows. First we utilise the triangle inequality to obtain,
    \begin{equation}
        W(\eta_{h},\eta_{h/2})\leq W(\eta_{h},\eta_{\star})+W(\eta_{\star},\eta_{h/2})
    \end{equation}
    Together with \cref{eq:conv_rate_assump_impl} this yields,
    \begin{equation}
        W(\eta_{h},\eta_{h/2})=Ch^{p}+C(h/2)^{p}+\mathcal{O}(h^{p+1})
    \end{equation}
    Similarly we have,
    \begin{equation}
        W(\eta_{h/2},\eta_{h/4})=C(h/2)^{p}+C(h/4)^{p}+\mathcal{O}(h^{p+1})
    \end{equation}
    Dividing the two above equations yields,
    \begin{align}
        \frac{W(\eta_{h},\eta_{h/2})}{W(\eta_{h/2},\eta_{h/4})} &= \frac{Ch^{p}+C(h/2)^{p}+\mathcal{O}(h^{p+1})}{C(h/2)^{p}+C(h/4)^{p}+\mathcal{O}(h^{p+1})} \nonumber \\
        &= \frac{1+2^{-p}+\mathcal{O}(h)}{2^{-p}+2^{-2p}+\mathcal{O}(h)} \nonumber \\
        &= 2^{p}+\mathcal{O}(h)
    \end{align}
    from which it follows that,
    \begin{equation}
        \operatorname{LR}(h):=\log_{2}\frac{W(\eta_{h},\eta_{h/2})}{W(\eta_{h/2},\eta_{h/4})}=p+\mathcal{O}(h)
    \end{equation}
    From the above equation we can see that we can obtain an estimate of the rate of convergence by computing, for a range of sufficiently small $h$-values, the base-2 logarithm of the ratio $W(\eta_{h},\eta_{h/2})/W(\eta_{h/2},\eta_{h/4})$ and seeing what these logarithms tend to as $h$ gets very small. \mylozenge
\end{remark}

As outlined in \cref{rem:conv_rate_true_res_unknown} above we will compute, for a range of sufficiently small $h$-values, the base-2 logarithm of the ratio $W(\nu_{h},\nu_{h/2})/W(\nu_{h/2},\nu_{h/4})$ and see what these logarithms tend to as $h$ approaches 0. From the results of \cref{thm:statfem_prior_error} we expect the logarithms to converge to 2.

\begin{figure}[!ht]
    \centering
    \includegraphics[width=12cm]{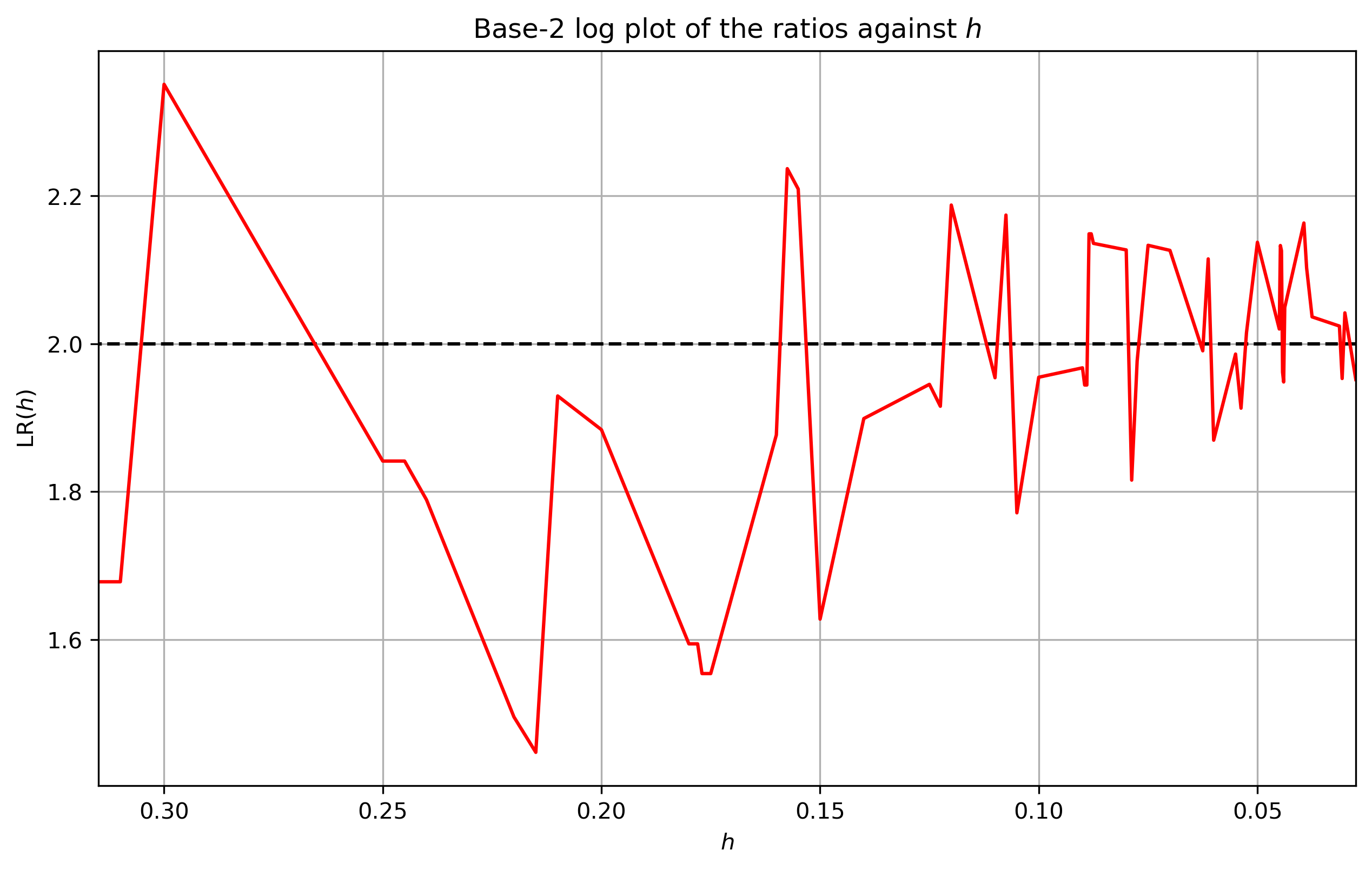}
    \caption{Plot of $\operatorname{LR}(h)$ against $h$, in the prior case, for our 2-D example.}
    \label{fig:2D_prior_results}
\end{figure}

We take a reference grid of $M=41^{2}=1681$ equally spaced points in $[0,1]^{2}$ to compare the covariances. We compute the logarithms for 59 $h$-values in the interval $[0.0275,0.315]$. Note that these $h$-values are not equally spaced; instead they are arranged according to a refinement strategy which was utilised for efficient use of memory. With these choices we obtain the plot shown in \cref{fig:2D_prior_results} on the next page.

\begin{figure}[!ht]
    \centering
    \includegraphics[width=12cm]{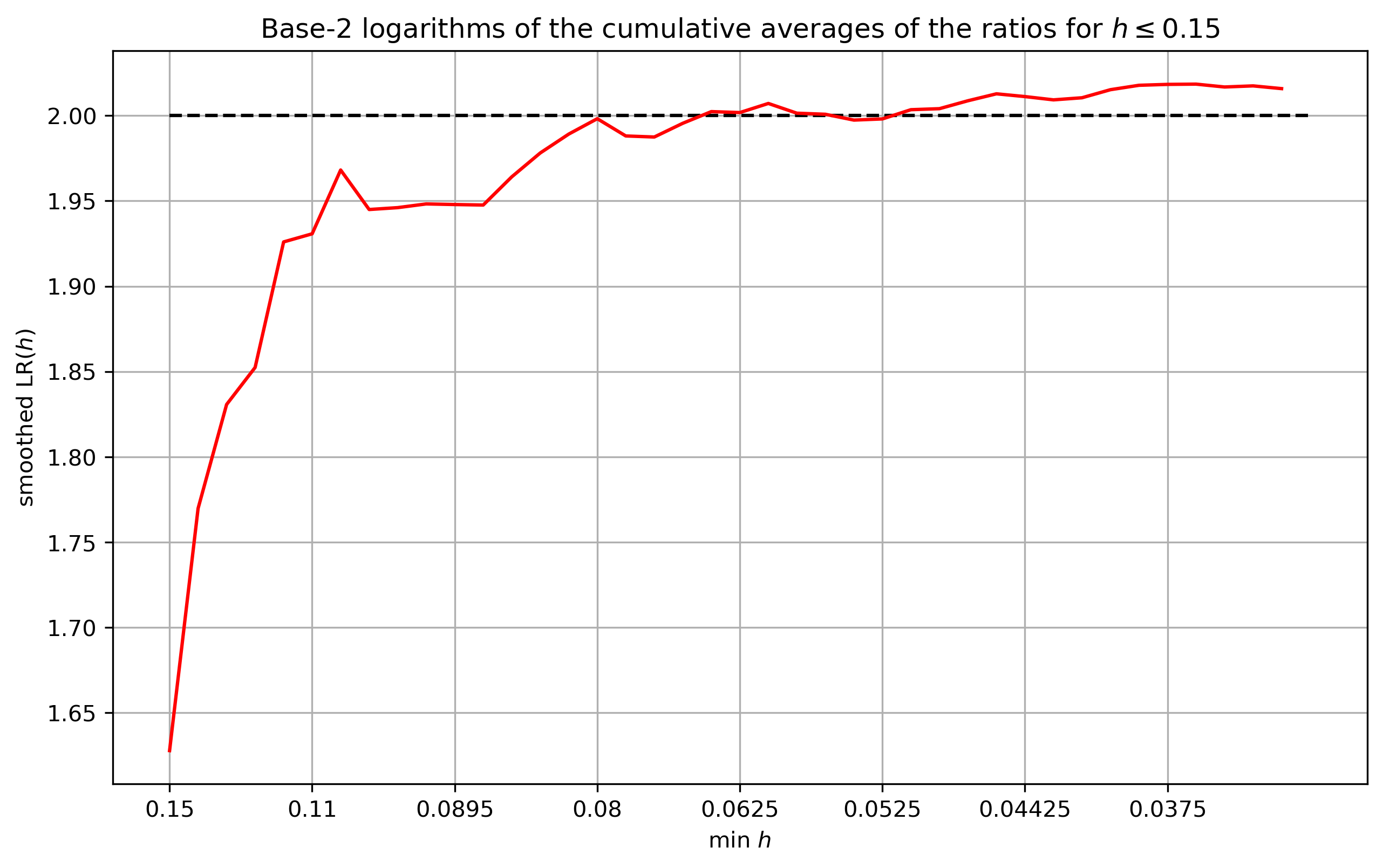}
    \caption{Plot of smoothed $\operatorname{LR}(h)$ for $h\leq 0.15$, in the prior case. The horizontal axis indicates the smallest $h$ value taken in the cumulative averages.}
    \label{fig:2D_prior_results_smoothed}
\end{figure}

From this plot we can see that the logarithms seem to be approaching $p=2$ as $h \to 0$, as our theory suggests. However, the results are somewhat oscillatory, which is to be expected for the approach outlined in \cref{rem:conv_rate_true_res_unknown}. Nevertheless the results seem to be oscillating around 2. To illustrate this more we smooth the above results by first discarding the results for larger $h$-values and then computing a cumulative average of the ratios before applying the base-2 logarithm. Note that this cumulative average is taken by starting with the results for large $h$ first. We take our cutoff point to be $h=0.15$. These smoothed results are shown in \cref{fig:2D_prior_results_smoothed}.

From \cref{fig:2D_prior_results_smoothed} we can see that by discarding the values corresponding to $h>0.15$ the base-2 logarithms of the rolling averages converge to a value slightly greater than $2$. The final logarithm for these rolling averages is $p=2.02$ (to 2 decimal places). Thus, we obtain results in close agreement with \cref{thm:statfem_prior_error}.

\subsubsection{Posterior results for 2-D example}\label{sec:2D_poisson_posterior}

For this example we will consider the same specific case of the general Bayesian inverse problem which we chose in \cref{sec:1D_poisson_posterior}. That is, we will again take the observation operator $S$ to be an evaluation operator. Since we do not have the true prior in an explicit form, we also do not have the true posterior and as such we will utilise the argument in \cref{rem:conv_rate_true_res_unknown} to demonstrate the \ac{statfem} posterior error bound.
To generate sensor readings, since we do not have access to the true prior, we instead simulated trajectories from the \ac{statfem} prior with a very fine grid of mesh size $h=1/120$.

We take a reference grid of $M=41^{2}=1681$ equally spaced points in $[0,1]^{2}$ to compare the covariances. We set the sensor noise level to $\epsilon=0.001$ and take $s=25$ sensor readings equally spaced in $[0.01,0.99]^{2}\subset D$. The level of sensor noise is chosen to be comparable to the prior variances. We follow the notation of \cref{thm:statfem_posterior_error} and denote the \ac{statfem} posterior by $\tilde{\nu}_{h}$. We compute the logarithms $\log_{2}\left(W(\tilde{\nu}_{h},\tilde{\nu}_{h/2})/W(\tilde{\nu}_{h/2},\tilde{\nu}_{h/4})\right)$ for 53 $h$-values in the interval $[0.03125,0.315]$, which are spread out according to the refinement strategy mentioned in \cref{sec:2D_poisson_prior}. With these choices we obtain the plot shown in \cref{fig:2D_posterior_results} below.

\begin{figure}[!ht]
    \centering
    \includegraphics[width=12cm]{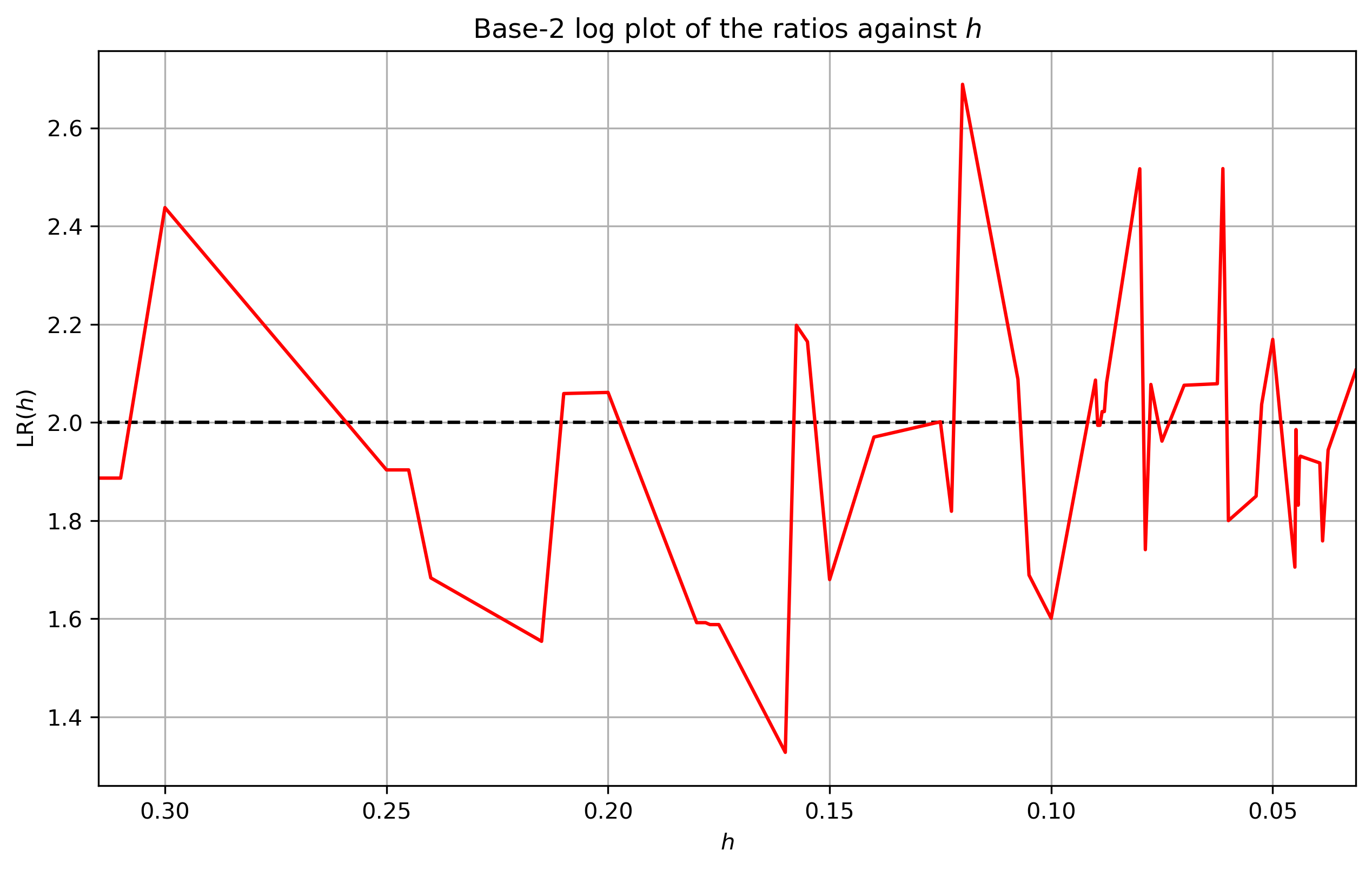}
    \caption{Plot of $\operatorname{LR}(h)$ against $h$, in the posterior case, for our 2-D example.}
    \label{fig:2D_posterior_results}
\end{figure}

\begin{figure}[!ht]
    \centering
    \includegraphics[width=12cm]{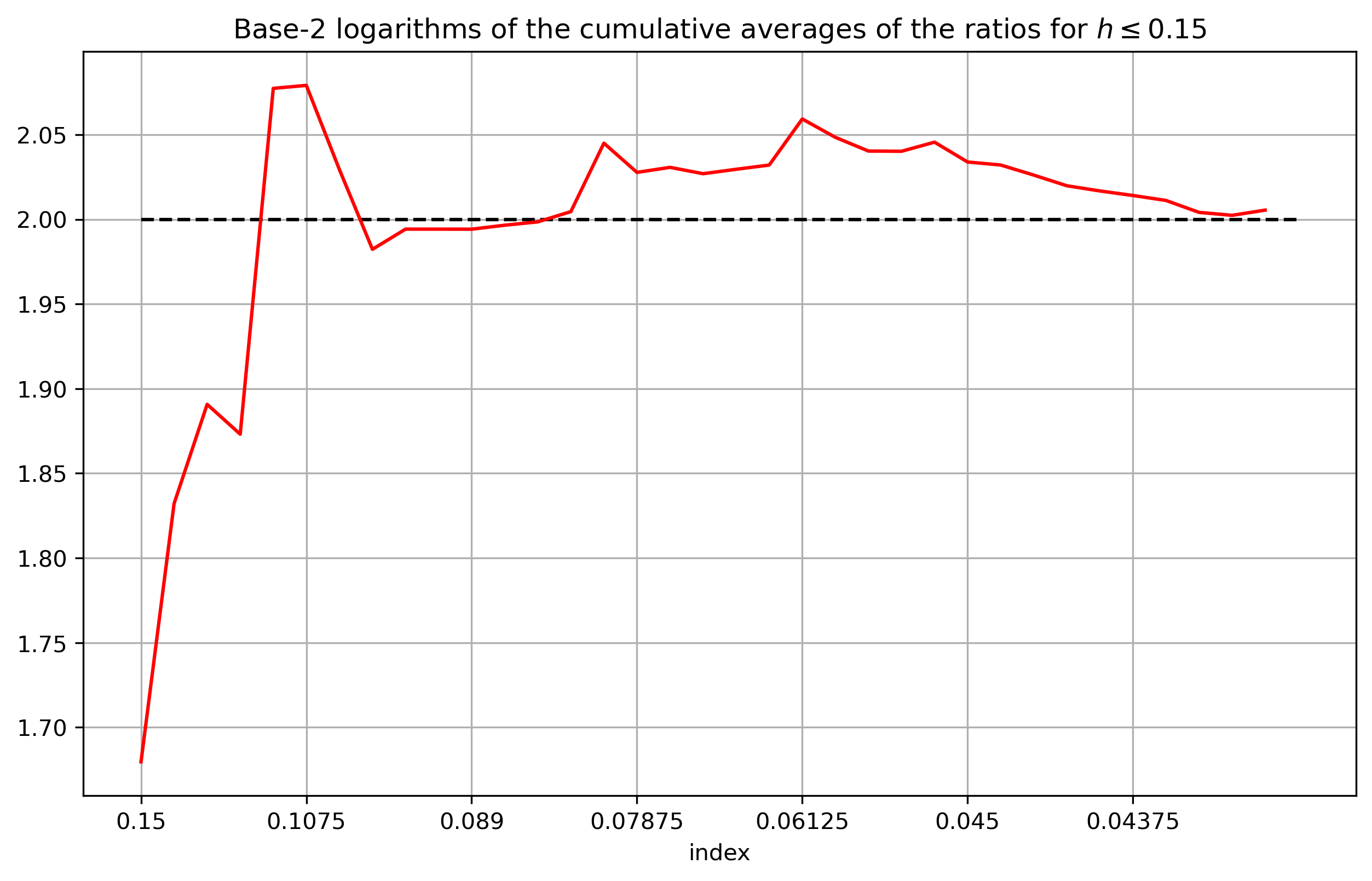}
    \caption{Plot of smoothed $\operatorname{LR}(h)$ for $h\leq 0.15$, in the posterior case. The horizontal axis indicates the smallest $h$ value taken in the cumulative averages.}
    \label{fig:2D_posterior_results_smoothed}
\end{figure}

From this plot we can see that the logarithms seem to be approaching $p=2$ as $h$ gets small. However, just as in \cref{sec:2D_poisson_prior}, the results are oscillatory but seem to be oscillating around $p=2$. We will thus smooth the results as explained in \cref{sec:2D_poisson_prior}. We again take the cutoff point to be $h=0.15$. These smoothed results are shown in \cref{fig:2D_posterior_results_smoothed} on the next page. From this figure we can see that by discarding the values corresponding to $h>0.15$ the base-2 logarithms of the rolling averages converge to a value slightly greater than 2. The final logarithm of these rolling averages is $p=2.00$ (to 2 decimal places). Thus, just as in the 1D example we see results which suggest that it is possible to extend \cref{thm:statfem_posterior_error} to the case of the full posteriors.

\subsection{Distribution of the maximum} \label{sec:1D_maximum}

In this section we will present a numerical experiment involving the same one-dimensional Poisson equation \cref{eq:1D_poisson_bvp}. The forcing $f$ will follow the same distribution as given in \cref{eq:1D_forcing_distr}, but we will adjust the kernel lengthscale to be $l_f=0.01$ instead, as with a larger length-scale the maximum of the \ac{statfem} prior and posterior were found to converge too quickly to demonstrate interesting behaviour. 
Since the solution, $u$, is random the maximum, $\max_{x\in[0,1]}u(x)$, is a univariate random variable. Since taking the maximum is not a linear functional we have no guarantees that the distribution of the maximum will be Gaussian. More importantly our results from \cref{thm:statfem_posterior_error} will not necessarily hold for this example. 

We will investigate the errors for the true and \ac{statfem} prior/posterior distributions for the maximum via a Monte Carlo approach as we now outline in the following remark. 
\begin{remark}
    \label{rem:max_dist_approach}
    Let $\eta_{\star},\eta_{h}$ denote the true and \ac{statfem} prior (or posterior) for the maximum of $u$ respectively. We are interested in investigating the rate of convergence of $W(\eta_{\star},\eta_{h})$ as $h\rightarrow 0$. This is done as follows: in both the prior and posterior case we will first fix a reference grid $\{\tilde{x_{i}}\}_{i=1}^{N}$ and obtain samples of trajectories of the solution on this grid from the true distribution. We then take the maximum values of these trajectories to obtain an (approximate) sample from the true distribution for $\max_{x\in[0,1]}u(x)$. For a range of $h$-values we will then simulate trajectories of the solution on the same reference grid from the \ac{statfem} distribution and then take the maximum values of these to obtain an (approximate) sample for the \ac{statfem} case. We will then compute an approximation of the Wasserstein distance $W(\eta_{\star},\eta_{h})$ by computing the Wasserstein distance between these two empirical distributions for each $h$-value by utilising the Python package \texttt{POT} \cite{flamary2021pot}. We then investigate the rate of convergence as outlined previously in \cref{rem:conv_rate_true_res_known}. \mylozenge
\end{remark}
 The results for the \ac{statfem} priors and posteriors will be presented separately in the following sections. 

\subsection{Prior results for maximum example}

As outlined in \cref{rem:max_dist_approach,rem:conv_rate_true_res_known}, we compute, for a range of sufficiently small $h$-values, a Monte Carlo approximation to the Wasserstein distance between the two priors for the maximum. We take a reference grid of $N=100$ equally spaced points in $[0,1]$ and $30$ $h$-values in $[0.02,0.25]$. We simulate $1000$ trajectories from both priors to obtain approximate samples of the maximum. With these choices we obtain the plot shown in \cref{fig:1D_max_prior_results} below. From this plot we can see that the we indeed have convergence of the Wasserstein distance for this non-linear quantity of interest to $0$, but at a different rate than that given in \cref{thm:statfem_prior_error}. The slope is estimated to be $p=1.35$ (to 2 decimal places).

\begin{figure}[!ht]
    \centering
    \includegraphics[width=12cm]{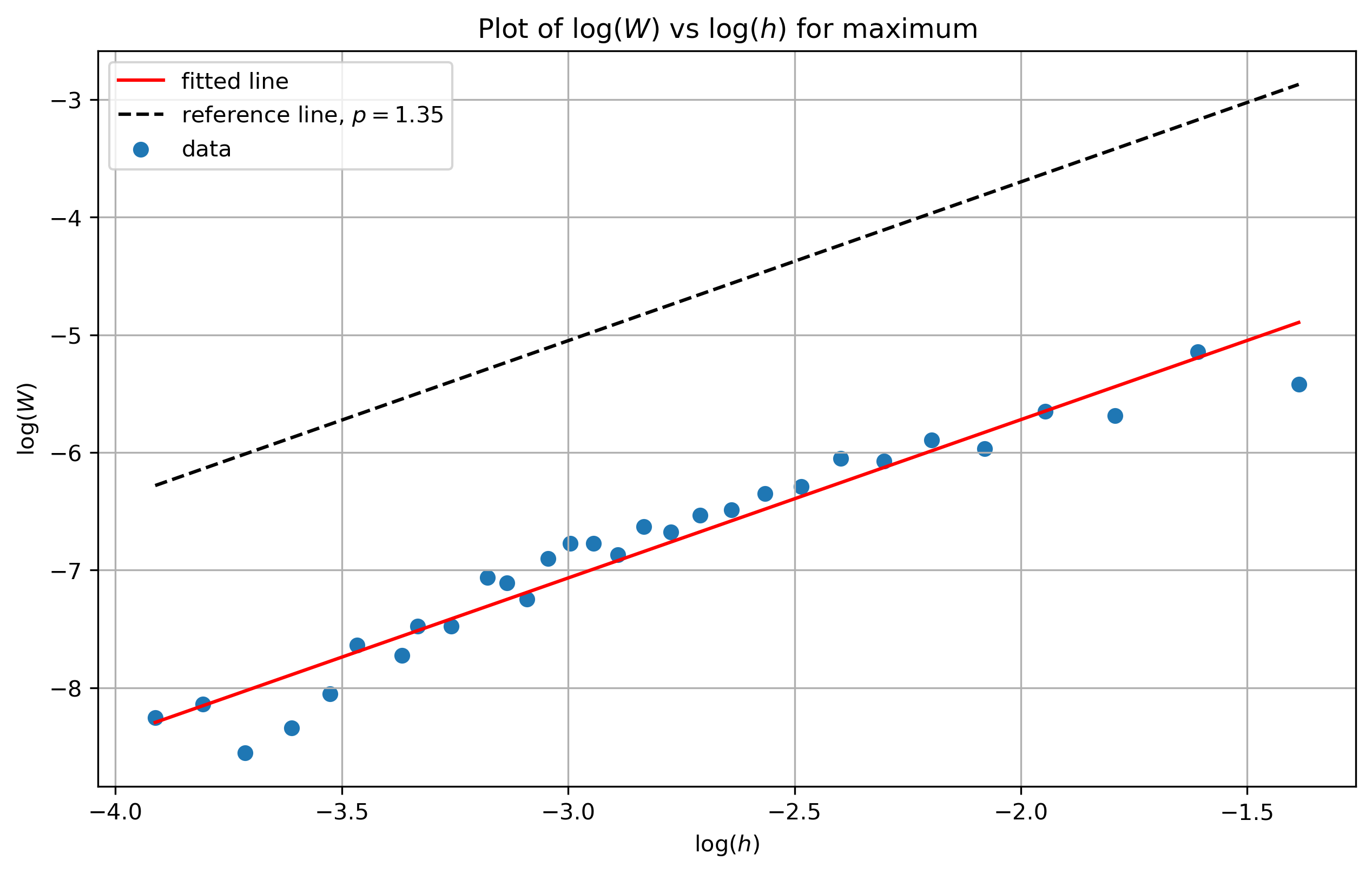}
    \caption{Log-log plot of Wasserstein distance against $h$ for the priors for the maximum. The reference line has a slope of 1.35. The slope of the fitted line is estimated to be $p=1.35$, in agreement visually with the reference line.}
    \label{fig:1D_max_prior_results}
\end{figure}

\subsection{Posterior results for maximum example}

\begin{figure}[!ht]
    \centering
    \includegraphics[width=12cm]{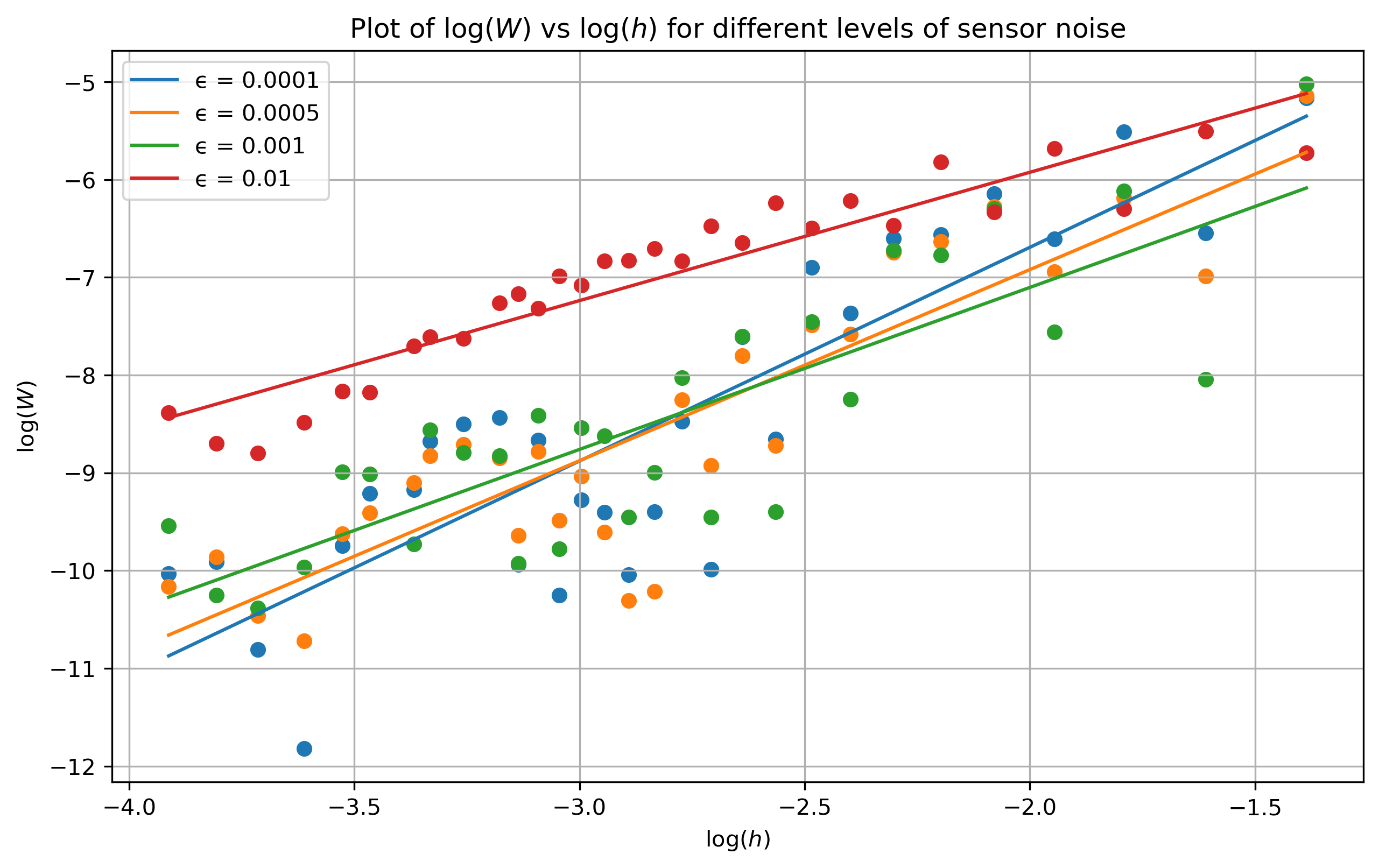}
    \caption{Log-log plot of Wasserstein distance against $h$ for the posteriors for the maximum for different levels of sensor noise.}
    \label{fig:1D_max_posterior_results}
\end{figure}

For the posterior case we take a reference grid of $N=41$ equally spaced points in $[0,1]$ to simulate $1000$ trajectories and $30$ $h$-values in $[0.02,0.25]$. We take $s=10$ sensor readings equally spaced in the interval $[0.01,0.99]$. This is repeated for $4$ different levels of sensor noise, $\epsilon=0.0001,0.0005,0.001,0.01$, again chosen so that the lower noise levels are comparable to the prior variances. For each level of noise we follow the argument outlined in \cref{rem:max_dist_approach,rem:conv_rate_true_res_known}. With these choices we obtain the plot shown in \cref{fig:1D_max_posterior_results} on the next page. We also present the estimates for the slopes and intercepts of the lines of best fit in \cref{tab:1D_max_posterior_results}.

\begin{table}[tbhp]
    {\footnotesize
    \caption{Estimate of slopes and intercepts for the different levels of sensor noise $\epsilon$.}\label{tab:1D_max_posterior_results}
    \begin{center}
        \begin{tabular}{@{}ccc@{}}
            \toprule
            $\epsilon$ & slope  & intercept \\ \midrule
            0.00010      & 2.1862 & -2.3203   \\
            0.00050       & 1.9552 & -3.0104   \\
            0.00100       & 1.6573 & -3.7885   \\
            0.01000       & 1.3131 & -3.2989   \\ \bottomrule
        \end{tabular}
    \end{center}
    }
\end{table}

From \cref{fig:1D_max_posterior_results} and \cref{tab:1D_max_posterior_results} we can see that for each level of sensor noise we do indeed have convergence to $0$ as $h$ goes to 0. However, we observe that the rate is not always the same as that given in \cref{thm:statfem_posterior_error}. From \cref{tab:1D_max_posterior_results} it is particularly clear that as the sensor noise level decreases the rate of convergence increases.
In particular, for $\epsilon=0.0001,0.0005$ we have rates around $p=2$, similar to the case of linear quantities of interest. We hypothesise that this is because as the sensor noise decreases the posterior distribution concentrates more as the solution is better identified and this, combined with the fact that we take sensor readings close to the where the true prior mean \cref{eq:1D_sol_mean} has a maximum, forces the maximum of the simulated trajectories to always be at the argmax of the prior mean in the exact case, i.e.\ at $x=0.5$. Thus for very low levels of sensor noise taking the maximum of the posterior draws is almost identical to evaluating the solution at $x=0.5$; such an evaluation is a linear quantity of interest and so one would expect the results from \cref{thm:statfem_posterior_error} to hold.

\section{Conclusion} 
\label{sec:conclusion}
In this paper we have provided detailed error analysis for the \ac{statfem} method first introduced in \cite{girolami2019statistical,Girolami2021} for the case of spatial elliptic \acp{bvp}. In particular, we have proven that the \ac{statfem} prior converges, with respect to the Wasserstein-2 distance, to the true prior at a rate of $\mathcal{O}(h^{2})$. We also showed that in the posterior case linear quantities of interest are also controlled at a rate of $\mathcal{O}(h^{2})$. In both settings, the choice of the Wasserstein distance provides a convenient metric in which to quantify the error introduced through discretisation. Numerical experiments demonstrating the theoretical results on standard \ac{fem} test problems were then presented. These demonstrated the rates for the prior theory, while also showcasing that it may be possible to extend the posterior results to the full distribution as opposed to just linear quantities of interest. An example involving a non-linear quantity of interest was also presented, showcasing an example which lies out of the remit of our theory.

Several possible avenues for future work present themselves. The first would be to prove error bounds for the full posterior distributions, extending \cref{thm:statfem_posterior_error}. Another avenue would be to provide similar error analysis for \ac{statfem} in the case of time-dependent problems. It would also be interesting to see if similar error bounds can be proven for distributions utilising approximations other than \ac{fem}. Generalising the results to utilise Wasserstein-$p$ distance, for $p>2$, would also be an interesting goal, as this will allow one to have even more control over other aspects of the distribution such as quantiles and extrema, which is especially important for risk-assessment applications.

\bibliographystyle{abbrvnat}
\bibliography{references}

\begin{thebibliography}{35}
\providecommand{\natexlab}[1]{#1}
\providecommand{\url}[1]{\texttt{#1}}
\expandafter\ifx\csname urlstyle\endcsname\relax
  \providecommand{\doi}[1]{doi: #1}\else
  \providecommand{\doi}{doi: \begingroup \urlstyle{rm}\Url}\fi

\bibitem[con()]{convergenceRateSite}
Verifying numerical convergence rates.
\newblock Technical report, KTH.
\newblock
  \url{https://www.csc.kth.se/utbildning/kth/kurser/DN2255/ndiff13/ConvRate.pdf}.

\bibitem[Aldosary et~al.(2018)Aldosary, Wang, and Li]{aldosary2018structural}
M.~Aldosary, J.~Wang, and C.~Li.
\newblock Structural reliability and stochastic finite element methods.
\newblock \emph{Engineering Computations}, 2018.

\bibitem[Bogachev(1998)]{bogachev1998gaussian}
V.~I. Bogachev.
\newblock \emph{Gaussian measures}.
\newblock Number~62. American Mathematical Soc., 1998.

\bibitem[Chkrebtii et~al.(2016)Chkrebtii, Campbell, Calderhead, Girolami,
  et~al.]{chkrebtii2016bayesian}
O.~A. Chkrebtii, D.~A. Campbell, B.~Calderhead, M.~A. Girolami, et~al.
\newblock Bayesian solution uncertainty quantification for differential
  equations.
\newblock \emph{Bayesian Analysis}, 11\penalty0 (4):\penalty0 1239--1267, 2016.

\bibitem[Cockayne et~al.(2016)Cockayne, Oates, Sullivan, and
  Girolami]{cockayne2016probabilistic}
J.~Cockayne, C.~Oates, T.~Sullivan, and M.~Girolami.
\newblock Probabilistic numerical methods for partial differential equations
  and bayesian inverse problems.
\newblock \emph{arXiv preprint arXiv:1605.07811}, 2016.

\bibitem[Cockayne et~al.(2019)Cockayne, Oates, Sullivan, and
  Girolami]{cockayne2019bayesian}
J.~Cockayne, C.~J. Oates, T.~J. Sullivan, and M.~Girolami.
\newblock Bayesian probabilistic numerical methods.
\newblock \emph{SIAM Review}, 61\penalty0 (4):\penalty0 756--789, 2019.

\bibitem[Conrad et~al.(2017)Conrad, Girolami, S{\"a}rkk{\"a}, Stuart, and
  Zygalakis]{conrad2017statistical}
P.~R. Conrad, M.~Girolami, S.~S{\"a}rkk{\"a}, A.~Stuart, and K.~Zygalakis.
\newblock Statistical analysis of differential equations: introducing
  probability measures on numerical solutions.
\newblock \emph{Statistics and Computing}, 27\penalty0 (4):\penalty0
  1065--1082, 2017.

\bibitem[Cotter et~al.(2009)Cotter, Dashti, Robinson, and
  Stuart]{cotter2009bayesian}
S.~L. Cotter, M.~Dashti, J.~C. Robinson, and A.~M. Stuart.
\newblock Bayesian inverse problems for functions and applications to fluid
  mechanics.
\newblock \emph{Inverse problems}, 25\penalty0 (11):\penalty0 115008, 2009.

\bibitem[Cuesta-Albertos et~al.(1996)Cuesta-Albertos, Matr{\'a}n-Bea, and
  Tuero-Diaz]{cuesta1996lower}
J.~Cuesta-Albertos, C.~Matr{\'a}n-Bea, and A.~Tuero-Diaz.
\newblock On lower bounds for the l2-wasserstein metric in a hilbert space.
\newblock \emph{Journal of Theoretical Probability}, 9\penalty0 (2):\penalty0
  263--283, 1996.

\bibitem[Da~Prato and Zabczyk(2014)]{da2014stochastic}
G.~Da~Prato and J.~Zabczyk.
\newblock \emph{Stochastic equations in infinite dimensions}.
\newblock Cambridge university press, 2014.

\bibitem[Dowson and Landau(1982)]{dowson1982frechet}
D.~Dowson and B.~Landau.
\newblock The fr{\'e}chet distance between multivariate normal distributions.
\newblock \emph{Journal of multivariate analysis}, 12\penalty0 (3):\penalty0
  450--455, 1982.

\bibitem[Evans(2010)]{evans10}
L.~C. Evans.
\newblock \emph{Partial differential equations}.
\newblock American Mathematical Society, Providence, R.I., 2010.
\newblock ISBN 9780821849743 0821849743.

\bibitem[Flamary et~al.(2021)Flamary, Courty, Gramfort, Alaya, Boisbunon,
  Chambon, Chapel, Corenflos, Fatras, Fournier, Gautheron, Gayraud, Janati,
  Rakotomamonjy, Redko, Rolet, Schutz, Seguy, Sutherland, Tavenard, Tong, and
  Vayer]{flamary2021pot}
R.~Flamary, N.~Courty, A.~Gramfort, M.~Z. Alaya, A.~Boisbunon, S.~Chambon,
  L.~Chapel, A.~Corenflos, K.~Fatras, N.~Fournier, L.~Gautheron, N.~T. Gayraud,
  H.~Janati, A.~Rakotomamonjy, I.~Redko, A.~Rolet, A.~Schutz, V.~Seguy, D.~J.
  Sutherland, R.~Tavenard, A.~Tong, and T.~Vayer.
\newblock Pot: Python optimal transport.
\newblock \emph{Journal of Machine Learning Research}, 22\penalty0
  (78):\penalty0 1--8, 2021.
\newblock URL \url{http://jmlr.org/papers/v22/20-451.html}.

\bibitem[Girolami et~al.(2019)Girolami, Gregory, Yin, and
  Cirak]{girolami2019statistical}
M.~Girolami, A.~Gregory, G.~Yin, and F.~Cirak.
\newblock The statistical finite element method.
\newblock \emph{arXiv preprint arXiv:1905.06391}, 2019.

\bibitem[Girolami et~al.(2021)Girolami, Febrianto, Yin, and
  Cirak]{Girolami2021}
M.~Girolami, E.~Febrianto, G.~Yin, and F.~Cirak.
\newblock The statistical finite element method ({statFEM}) for coherent
  synthesis of observation data and model predictions.
\newblock \emph{Computer Methods in Applied Mechanics and Engineering},
  375:\penalty0 113533, Mar. 2021.
\newblock \doi{10.1016/j.cma.2020.113533}.
\newblock URL \url{https://doi.org/10.1016/j.cma.2020.113533}.

\bibitem[Grisvard(2011)]{grisvard2011elliptic}
P.~Grisvard.
\newblock \emph{Elliptic problems in nonsmooth domains}.
\newblock SIAM, 2011.

\bibitem[Hennig et~al.(2015)Hennig, Osborne, and
  Girolami]{hennig2015probabilistic}
P.~Hennig, M.~A. Osborne, and M.~Girolami.
\newblock Probabilistic numerics and uncertainty in computations.
\newblock \emph{Proceedings of the Royal Society A: Mathematical, Physical and
  Engineering Sciences}, 471\penalty0 (2179):\penalty0 20150142, 2015.

\bibitem[Larsson and Thom{\'e}e(2008)]{larsson2008partial}
S.~Larsson and V.~Thom{\'e}e.
\newblock \emph{Partial differential equations with numerical methods},
  volume~45.
\newblock Springer Science \& Business Media, 2008.

\bibitem[Lindgren et~al.(2011)Lindgren, Rue, and
  Lindstr{\"o}m]{lindgren2011explicit}
F.~Lindgren, H.~Rue, and J.~Lindstr{\"o}m.
\newblock An explicit link between gaussian fields and gaussian markov random
  fields: the stochastic partial differential equation approach.
\newblock \emph{Journal of the Royal Statistical Society: Series B (Statistical
  Methodology)}, 73\penalty0 (4):\penalty0 423--498, 2011.

\bibitem[Lord et~al.(2014)Lord, Powell, and Shardlow]{lord2014introduction}
G.~J. Lord, C.~E. Powell, and T.~Shardlow.
\newblock \emph{An introduction to computational stochastic PDEs}, volume~50.
\newblock Cambridge University Press, 2014.

\bibitem[Martino and Riebler(2014)]{martino2014integrated}
S.~Martino and A.~Riebler.
\newblock Integrated nested laplace approximations (inla).
\newblock \emph{Wiley StatsRef: Statistics Reference Online}, pages 1--19,
  2014.

\bibitem[Masarotto et~al.(2019)Masarotto, Panaretos, and
  Zemel]{masarotto2019procrustes}
V.~Masarotto, V.~M. Panaretos, and Y.~Zemel.
\newblock Procrustes metrics on covariance operators and optimal transportation
  of gaussian processes.
\newblock \emph{Sankhya A}, 81\penalty0 (1):\penalty0 172--213, 2019.

\bibitem[Matthies et~al.(1997)Matthies, Brenner, Bucher, and {Guedes
  Soares}]{MATTHIES1997283}
H.~G. Matthies, C.~E. Brenner, C.~G. Bucher, and C.~{Guedes Soares}.
\newblock Uncertainties in probabilistic numerical analysis of structures and
  solids-stochastic finite elements.
\newblock \emph{Structural Safety}, 19\penalty0 (3):\penalty0 283--336, 1997.
\newblock ISSN 0167-4730.
\newblock \doi{https://doi.org/10.1016/S0167-4730(97)00013-1}.
\newblock URL
  \url{https://www.sciencedirect.com/science/article/pii/S0167473097000131}.
\newblock Devoted to the work of the Joint Committee on Structural Safety.

\bibitem[Olkin and Pukelsheim(1982)]{olkin1982distance}
I.~Olkin and F.~Pukelsheim.
\newblock The distance between two random vectors with given dispersion
  matrices.
\newblock \emph{Linear Algebra and its Applications}, 48:\penalty0 257--263,
  1982.

\bibitem[Owhadi(2015)]{owhadi2015bayesian}
H.~Owhadi.
\newblock Bayesian numerical homogenization.
\newblock \emph{Multiscale Modeling \& Simulation}, 13\penalty0 (3):\penalty0
  812--828, 2015.

\bibitem[Rasmussen(2006)]{Rasmussen06gaussianprocesses}
C.~E. Rasmussen.
\newblock \emph{Gaussian processes for machine learning}.
\newblock MIT Press, 2006.

\bibitem[Scheuerer(2010)]{scheuerer2010regularity}
M.~Scheuerer.
\newblock Regularity of the sample paths of a general second order random
  field.
\newblock \emph{Stochastic Processes and their Applications}, 120\penalty0
  (10):\penalty0 1879--1897, 2010.

\bibitem[Sigrist et~al.(2015)Sigrist, K{\"u}nsch, and
  Stahel]{sigrist2015stochastic}
F.~Sigrist, H.~R. K{\"u}nsch, and W.~A. Stahel.
\newblock Stochastic partial differential equation based modelling of large
  space--time data sets.
\newblock \emph{Journal of the Royal Statistical Society: Series B: Statistical
  Methodology}, pages 3--33, 2015.

\bibitem[Stefanou(2009)]{STEFANOU20091031}
G.~Stefanou.
\newblock The stochastic finite element method: Past, present and future.
\newblock \emph{Computer Methods in Applied Mechanics and Engineering},
  198\penalty0 (9):\penalty0 1031--1051, 2009.
\newblock ISSN 0045-7825.
\newblock \doi{https://doi.org/10.1016/j.cma.2008.11.007}.
\newblock URL
  \url{https://www.sciencedirect.com/science/article/pii/S0045782508004118}.

\bibitem[Strang and Fix(1973)]{strang1973analysis}
G.~Strang and G.~J. Fix.
\newblock An analysis of the finite element method.
\newblock \emph{Journal of Applied Mathematics and Mechanics}, 1973.

\bibitem[Stuart(2010)]{stuart2010inverse}
A.~M. Stuart.
\newblock Inverse problems: a bayesian perspective.
\newblock \emph{Acta numerica}, 19:\penalty0 451--559, 2010.

\bibitem[Sudret and Der~Kiureghian(2000)]{sudret2000stochastic}
B.~Sudret and A.~Der~Kiureghian.
\newblock \emph{Stochastic finite element methods and reliability: a
  state-of-the-art report}.
\newblock Department of Civil and Environmental Engineering, University of
  California~…, 2000.

\bibitem[Tronarp et~al.(2019)Tronarp, Kersting, S{\"a}rkk{\"a}, and
  Hennig]{tronarp2019probabilistic}
F.~Tronarp, H.~Kersting, S.~S{\"a}rkk{\"a}, and P.~Hennig.
\newblock Probabilistic solutions to ordinary differential equations as
  nonlinear bayesian filtering: a new perspective.
\newblock \emph{Statistics and Computing}, 29\penalty0 (6):\penalty0
  1297--1315, 2019.

\bibitem[Xiu(2010)]{xiu2010numerical}
D.~Xiu.
\newblock \emph{Numerical methods for stochastic computations: a spectral
  method approach}.
\newblock Princeton university press, 2010.

\bibitem[Zhang and Karniadakis(2017)]{zhang2017numerical}
Z.~Zhang and G.~Karniadakis.
\newblock \emph{Numerical methods for stochastic partial differential equations
  with white noise}.
\newblock Springer, 2017.

\end{thebibliography}


\begin{thebibliography}{4}
\providecommand{\natexlab}[1]{#1}
\providecommand{\url}[1]{\texttt{#1}}
\expandafter\ifx\csname urlstyle\endcsname\relax
  \providecommand{\doi}[1]{doi: #1}\else
  \providecommand{\doi}{doi: \begingroup \urlstyle{rm}\Url}\fi

\bibitem[Cockayne and Duncan(2020)]{cockayne2020probabilistic}
J.~Cockayne and A.~B. Duncan.
\newblock Probabilistic gradients for fast calibration of differential equation
  models.
\newblock \emph{arXiv preprint arXiv:2009.04239}, 2020.

\bibitem[Girolami et~al.(2021)Girolami, Febrianto, Yin, and
  Cirak]{Girolami2021}
M.~Girolami, E.~Febrianto, G.~Yin, and F.~Cirak.
\newblock The statistical finite element method ({statFEM}) for coherent
  synthesis of observation data and model predictions.
\newblock \emph{Computer Methods in Applied Mechanics and Engineering},
  375:\penalty0 113533, Mar. 2021.
\newblock \doi{10.1016/j.cma.2020.113533}.
\newblock URL \url{https://doi.org/10.1016/j.cma.2020.113533}.

\bibitem[Gohberg et~al.(2012)Gohberg, Goldberg, and Kaashoek]{gohberg2012basic}
I.~Gohberg, S.~Goldberg, and M.~Kaashoek.
\newblock \emph{Basic classes of linear operators}.
\newblock Birkh{\"a}user, 2012.

\bibitem[Stuart(2010)]{stuart2010inverse}
A.~M. Stuart.
\newblock Inverse problems: a bayesian perspective.
\newblock \emph{Acta numerica}, 19:\penalty0 451--559, 2010.

\end{thebibliography}



\section*{Supplementary material}

# Supplementary Material

## For: Theoretical Guarantees for the Statistical Finite Element Method

Yanni Papandreou[*]   Jon Cockayne [†]   Mark Girolami [‡]

Andrew B. Duncan[§]

February 18, 2022

## S1 Implementation Details for Experiments

In this section we present three remarks which give further details on the implementation of the experiments in Section 4.

***Remark*** 1. In both the statistical finite element method (StatFEM) prior and posterior the $n_u \times n_u$ matrix $K_{\mathcal{F}_h} = \mathcal{F}_h K \mathcal{F}_h^\dagger$ appears. This is the covariance matrix for the random load vector $\mathbf{f} = \mathcal{F}_h f$. The $ij$-*th* element of this covariance matrix is given by:

$$(K_{\mathcal{F}_h})_{ij} = \int_D \int_D \phi_i(x) k_f(x,y) \phi_j(y) \mathrm{d}x \mathrm{d}y. \tag{S1}$$

For a covariance function $k_f$ which is globally supported, the covariance matrix, $K_{\mathcal{F}_h}$, is dense as can be seen from expression Eq. (S1) above. As such, it can be very difficult to efficiently compute and assemble the entries of this matrix. We thus employ an approximation of this matrix, by interpolating the covariance function $k_f$ with the finite element basis functions, as suggested by [S2]. I.e., denoting the non-boundary nodes of the FE triangulation by $\{x_i\}_{i=1}^{n_u}$, we approximate $k_f$ by,

$$k_f(x,y) \approx \sum_{k,l=1}^{n_u} \phi_k(x) k_f(x_k, x_l) \phi_l(y) \tag{S2}$$

---

[*]Imperial College London (john.papandreou18@imperial.ac.uk)

[†]The Alan Turing Institute (jcockayne@turing.ac.uk)

[‡]The Alan Turing Institute and University of Cambridge (mag92@cam.ac.uk)

[§]Imperial College London and The Alan Turing Institute (a.duncan@imperial.ac.uk)

Using this interpolation of $k_f$ we can approximate the entries of $K_{\mathcal{F}_h}$ by,

$$(K_{\mathcal{F}_h})_{ij} \approx \sum_{k,l=1}^{n_u} k_f(x_k, x_l) \int_D \int_D \phi_i(x)\phi_k(x)\phi_j(y)\phi_l(y) \mathrm{d}x\mathrm{d}y$$

$$= \sum_{k,l=1}^{n_u} k_f(x_k, x_l) \int_D \phi_i(x)\phi_k(x)\mathrm{d}x \int_D \phi_j(y)\phi_l(y)\mathrm{d}y$$

$$= \sum_{k,l=1}^{n_u} k_f(x_k, x_l) M_{ik} M_{jl}$$

$$= (MC_f M^T)_{ij}$$

i.e. we can approximate:
$$K_{\mathcal{F}_h} \approx MC_f M^T \tag{S3}$$

where $C_f \in \mathbb{R}^{n_u \times n_u}$ is the covariance matrix corresponding to evaluating $k_f$ on the non-boundary nodes, $\{x_i\}_{i=1}^{n_u}$, of the Finite Element Method (FEM) triangulation, and $M$ is the mass matrix introduced in Section 2.1. This approximation will improve as $h$ gets smaller. ♦

***Remark*** 2. Since all of the true and StatFEM priors and posteriors are Gaussian measures, and since we are interested in demonstrating a bound on the error between the respective priors and posteriors for the Wasserstein-2 distance the expression Eq. (3) will be useful. Focussing on the covariance contribution in Eq. (3), note that this is difficult to compute explicitly for the measures in this paper. For two Gaussian processes (GPs) $\nu_i \sim \mathcal{N}(0, \Sigma_i)$, $i \in \{1, 2\}$ we approximate the covariance contribution by fixing a fine reference grid $\{\tilde{x}_i\}_{i=1}^N$, computing the covariance matrices, $C_i = \{k_i(\tilde{x}_k, \tilde{x}_l)\}_{k,l=1}^N$, corresponding to the covariance operators, $\Sigma_i$, evaluated on this reference grid, where $k_i$ is the covariance function associated with $\Sigma_i$ respectively, and finally approximating the Wasserstein distance as follows:

$$W^2(\mathcal{N}(0, \Sigma_1), \mathcal{N}(0, \Sigma_2)) = \mathrm{tr}(\Sigma_1) + \mathrm{tr}(\Sigma_2) - 2\,\mathrm{tr}\sqrt{\Sigma_1^{1/2} \Sigma_2 \Sigma_1^{1/2}}$$
$$\approx \mathrm{tr}(C_1) + \mathrm{tr}(C_2) - 2\,\mathrm{tr}\sqrt{C_1^{1/2} C_2 C_1^{1/2}}. \tag{S4}$$

♦

***Remark*** 3. Let $m_\star, m_h$ denote the means of the true and StatFEM prior (or posterior) respectively. The mean contribution to the Wasserstein-2 distance is $\|m_\star - m_h\|_{L^2(D)}$. This contribution can be easily approximated by utilising a modern FEM package such as `FEniCS`, which we use in this paper. ♦

## S2 Proofs

In this section the proofs of the main theoretical results are presented.

*Proof of Proposition 1.* This follows trivially from the fact that the approximate solution operator is a linear map together with the well known transformation results of a Gaussian measure by linear maps. ∎

*Proof of Proposition 2.* We first note that so far we have assumed that the priors $u \sim \nu_i = \mathcal{N}(m_i, \Sigma_i)$ are Gaussian measures on $L^2(D)$. However, both priors are actually supported on subspaces $V_i \subseteq L^2(D)$. For the purposes of the proof we will be viewing the priors as measures on the respective subspaces $V_i$. We will then trivially extend the resulting posteriors onto all of $L^2(D)$. When viewed as Gaussian measures on $V_i$ the covariance operators $\Sigma_i$ are both positive definite.

In order to obtain the posterior $\widetilde{\nu}_i$ we will first derive the joint distribution of $(u, \mathbf{v})^T$ when the prior on $u$ is $\nu_i$. Since the observational noise $\boldsymbol{\xi}$ is independent of $\mathbf{u}$ we have that the joint distribution of $(u, \boldsymbol{\xi})^T$ is:

$$\begin{pmatrix} u \\ \boldsymbol{\xi} \end{pmatrix} \sim \mathcal{N}\left( \begin{pmatrix} m_i \\ \mathbf{0} \end{pmatrix}, \begin{pmatrix} \Sigma_i & 0 \\ 0 & \epsilon^2 I \end{pmatrix} \right). \tag{S5}$$

Using the definition Eq. (20) of $\mathbf{v}$ we now note that $(u, \mathbf{v})^T$ can be expressed as:

$$\begin{pmatrix} u \\ \mathbf{v} \end{pmatrix} = \begin{pmatrix} I & 0 \\ S & I \end{pmatrix} \begin{pmatrix} u \\ \boldsymbol{\xi} \end{pmatrix} \tag{S6}$$

Since this is a linear transformation of $(u, \boldsymbol{\xi})^T$ we have that the joint distribution of $(u, \mathbf{v})^T$ is given by:

$$\begin{pmatrix} u \\ \mathbf{v} \end{pmatrix} \sim \mathcal{N}\left( \begin{pmatrix} m_i \\ Sm_i \end{pmatrix}, \begin{pmatrix} \Sigma_i & \Sigma_i S^\dagger \\ S\Sigma_i & \epsilon^2 I + S\Sigma_i S^\dagger \end{pmatrix} \right) \tag{S7}$$

From this it is simply a matter of conditioning to obtain the posterior distribution of $u|\mathbf{v}$ specified by Eqs. (23) and (24). This conditioning step is justified under the following result taken from [S4, see Theorem 6.20].

**Theorem 1** (Conditioning of Gaussian measures)**.** *Let $\mathcal{H} = \mathcal{H}_1 \oplus \mathcal{H}_2$ be a separable Hilbert space with projectors $\Pi_i : \mathcal{H} \to \mathcal{H}_i$. Let $(x_1, x_2) \in \mathcal{H}_1 \oplus \mathcal{H}_2$ be an $\mathcal{H}$-valued Gaussian random variable with mean $m = (m_1, m_2)$ and positive definite covariance operator $\mathcal{C}$. Define*

$$C_{ij} = \mathbb{E}(x_i - m_i) \otimes (x_j - m_j) \tag{S8}$$

*Then the conditional distribution of $x_1$ given $x_2$ is Gaussian with mean*

$$m' = m_1 + \mathcal{C}_{12}\mathcal{C}_{22}^{-1}(x_2 - m_2) \tag{S9}$$

*and covariance operator*

$$\mathcal{C}' = \mathcal{C}_{11} - \mathcal{C}_{12}C_{22}^{-1}\mathcal{C}_{21} \tag{S10}$$

♦

This theorem is applied in our case with $\mathcal{H}_1 = V_i$ and $\mathcal{H}_2 = \mathbb{R}^s$. The corresponding covariance operator $\mathcal{C}$ is given by the operator in Eq. (S7). The only thing which remains to show is that this covariance operator is positive definite. To that end, note that we can rewrite this operator as:

$$\begin{aligned}
\mathcal{C} &= \begin{pmatrix} \Sigma_i & \Sigma_i S^\dagger \\ S\Sigma_i & \epsilon^2 I + S\Sigma_i S^\dagger \end{pmatrix} = \begin{pmatrix} \Sigma_i & \Sigma_i S^\dagger \\ S\Sigma_i & S\Sigma_i S^\dagger \end{pmatrix} + \begin{pmatrix} 0 & 0 \\ 0 & \epsilon^2 I \end{pmatrix} \\
&= \underbrace{\begin{pmatrix} I \\ S \end{pmatrix} \Sigma_i \begin{pmatrix} I & S^\dagger \end{pmatrix}}_{(a)} + \underbrace{\begin{pmatrix} 0 & 0 \\ 0 & \epsilon^2 I \end{pmatrix}}_{(b)}
\end{aligned} \quad (S11)$$

Written in this form it is clear that $\mathcal{C}$ is positive definite. This is because operator $(a)$ above is positive definite while operator $(b)$ is positive semi-definite. Positive definiteness of $(a)$ follows from the positive definiteness of $\Sigma_i$ together with the fact that $S$ satisfies Assumption 4. ∎

Before we provide proofs for Theorems 3 and 4 we first present several lemmas which will be required in the proofs. The first of these is a simple corollary to Theorem 2 which places a bound on the operator norm of the difference of the solution operators $\mathcal{L}^{-1}$ and $\hat{L}_h^{-1}$.

**Lemma 1.** Under the assumptions of Theorem 2 we have the following bound for sufficiently small $h$:
$$\|\mathcal{L}^{-1} - \hat{L}_h^{-1}\|_{L^2} \leq Ch^2 \quad (S12)$$
where we have used the notation $\|\cdot\|_{L^2}$ to denote the operator norm for operators from $L^2(D) \to L^2(D)$. ♦

*Proof of Lemma 1.* This bound on the operator norm of $\mathcal{L}^{-1} - \hat{L}_h^{-1}$ follows from noting that the error bound Eq. (14) given in Theorem 2 can be re-expressed as:
$$\|\mathcal{L}^{-1} f - \hat{L}_h^{-1} f\|_{L^2(D)} = \|(\mathcal{L}^{-1} - \hat{L}_h^{-1}) f\|_{L^2(D)} \leq Ch^2 \|f\|_{L^2(D)}$$

Since the above inequality holds for all $f \in L^2(D)$ we deduce that Eq. (S12) holds as required. ∎

The next result we will need is a bound on the Wasserstein distance found in [S1, see Lemma S2]. We quote this lemma below:

**Lemma 2.** Let $\mu_1 = \mathcal{N}(a_1, C_1)$ and $\mu_2 = \mathcal{N}(a_2, C_2)$ each be Gaussian measures on a Hilbert space $\mathcal{H}$. Let $S_1, S_2$ each be square roots of $C_1, C_2$ such that $C_1 = S_1^\dagger S_1$ and $C_2 = S_2^\dagger S_2$. Then we have that,
$$W^2(\mu_1, \mu_2) \leq \|a_1 - a_2\|_\mathcal{H}^2 + \|S_1 - S_2\|_{\text{HS}}^2 \quad (S13)$$

♦

4

***Remark*** 4. Lemma 2 above utilises the connection between the Wasserstein distance between Gaussian measures and the Procrustes Metric on the respective covariance operators mentioned above. ♦

We are now in a position to prove Theorems 3 and 4. We start by proving Theorem 3.

*Proof of Theorem 3.* Both of the priors $\nu_\star, \nu_h$ are Gaussian measures on the Hilbert space $L^2(D)$. Further, noting that $K$ is a covariance operator, we have that it possesses a symmetric square root $K^{1/2}$, i.e. we have $K = K^{1/2}K^{1/2}$ and $(K^{1/2})^\dagger = K^{1/2}$. Using this fact, we see that we can express the covariance operators of the priors as follows:

$$\begin{aligned}
\Sigma_\star &= \mathcal{L}^{-1} K \mathcal{L}^{-\dagger} \\
&= \mathcal{L}^{-1} K^{1/2} K^{1/2} \mathcal{L}^{-\dagger} \\
&= (K^{1/2}\mathcal{L}^{-\dagger})^\dagger (K^{1/2}\mathcal{L}^{-\dagger}) = S_\star^\dagger S_\star
\end{aligned} \tag{S14}$$

$$\begin{aligned}
\Sigma_h &= \hat{L}_h^{-1} K \hat{L}_h^{-\dagger} \\
&= \hat{L}_h^{-1} K^{1/2} K^{1/2} \hat{L}_h^{-\dagger} \\
&= (K^{1/2}\hat{L}_h^{-\dagger})^\dagger (K^{1/2}\hat{L}_h^{-\dagger}) = S_h^\dagger S_h
\end{aligned} \tag{S15}$$

where $S_\star := K^{1/2}\mathcal{L}^{-\dagger}$ and $S_h := K^{1/2}\hat{L}_h^{-\dagger}$. Thus, our priors $\nu_\star, \nu_h$ satisfy the assumptions of Lemma 2 and so utilising this lemma we obtain the following bound on the Wasserstein distance:

$$W^2(\nu_\star, \nu_h) \leq \|m_\star - m_h\|_{L^2(D)}^2 + \|S_\star - S_h\|_{\text{HS}}^2 \tag{S16}$$

We now proceed to bound each of the terms in Eq. (S16). Starting with the first term we have, for sufficiently small $h$, the following:

$$\|m_\star - m_h\|_{L^2(D)} = \|\mathcal{L}^{-1}\bar{f} - \hat{L}_h^{-1}\bar{f}\|_{L^2(D)} \leq Ch^2 \|\bar{f}\|_{L^2(D)} \tag{S17}$$

where we have utilised the error bound from Theorem 2, which holds since $\bar{f} \in L^2(D)$ and since all of the Assumptions 1 to 3 hold. For the second term in Eq. (S16) we have:

$$\|S_\star - S_h\|_{\text{HS}} = \|K^{1/2}\mathcal{L}^{-\dagger} - K^{1/2}\hat{L}_h^{-\dagger}\|_{\text{HS}} \tag{S18}$$

$$= \|K^{1/2}(\mathcal{L}^{-\dagger} - \hat{L}_h^{-\dagger})\|_{\text{HS}} \tag{S19}$$

$$\leq \|K^{1/2}\|_{HS} \|\mathcal{L}^{-\dagger} - \hat{L}_h^{-\dagger}\|_{L^2} \tag{S20}$$

$$= \|K^{1/2}\|_{\text{HS}} \|\mathcal{L}^{-1} - \hat{L}_h^{-1}\|_{L^2} \tag{S21}$$

$$\leq Ch^2 \|K^{1/2}\|_{\text{HS}} \tag{S22}$$

We have utilised several facts in the above which we now detail. In going from Eq. (S19) to Eq. (S20) we have utilised the facts that $\mathcal{L}^{-1} - \hat{L}_h^{-1}$, and hence its adjoint $\mathcal{L}^{-\dagger} - \hat{L}_h^{-\dagger}$, is a bounded operator and that $K^{1/2}$ is Hilbert-Schmidt. This second fact follows from the facts that $K$ and $K^{1/2}$ are symmetric and that $K$ is positive and trace-class and computing $\|K^{1/2}\|_{\text{HS}} = \text{tr}(K^{1/2}K^{1/2})^{1/2} = \text{tr}(K)^{1/2} = \text{tr}((K^2)^{1/2})^{1/2} = \|K\|_{\text{tr}}^{1/2} < \infty$. Finally, in going from Eq. (S21) to Eq. (S22) we used the bound Eq. (S12) from Lemma 1.

5

We can now combine the bounds on these two terms to obtain for sufficiently small $h$, the following:

$$W^2(\nu_\star, \nu_h) \leq C^2 h^4 \|\bar{f}\|^2_{L^2(D)} + C^2 h^4 \|K^{1/2}\|^2_{\text{HS}} \tag{S23}$$

$$= C^2 h^4 \left( \|\bar{f}\|^2_{L^2(D)} + \|K^{1/2}\|^2_{\text{HS}} \right) \tag{S24}$$

which implies,

$$W(\nu_\star, \nu_h) \leq C \sqrt{\|\bar{f}\|^2_{L^2(D)} + \|K^{1/2}\|^2_{\text{HS}}} h^2 = \gamma h^2 \tag{S25}$$

where $\gamma := C\sqrt{\|\bar{f}\|^2_{L^2(D)} + \|K^{1/2}\|^2_{\text{HS}}} > 0$ is a positive constant, independent of $h$, as required. $\blacksquare$

We now prove Theorem 4.

*Proof of Theorem 4.* Since both of the posteriors $\tilde{\nu}_\star, \tilde{\nu}_h$ are Gaussians we have that $\ell_\# \tilde{\nu}_\star \sim \mathcal{N}\left(\ell m^{(\star)}_{u|\mathbf{v}}, \ell \Sigma^{(\star)}_{u|\mathbf{v}} \ell^\dagger\right)$ and $\ell_\# \tilde{\nu}_h \sim \mathcal{N}\left(\ell m^{(h)}_{u|\mathbf{v}}, \ell \Sigma^{(h)}_{u|\mathbf{v}} \ell^\dagger\right)$ are univariate Gaussians. Thus, we have via Eq. (3),

$$W^2(\ell_\# \tilde{\nu}_\star, \ell_\# \tilde{\nu}_h) = \left| \ell m^{(\star)}_{u|\mathbf{v}} - \ell m^{(h)}_{u|\mathbf{v}} \right|^2 + \left| \sqrt{\ell \Sigma^{(\star)}_{u|\mathbf{v}} \ell^\dagger} - \sqrt{\ell \Sigma^{(h)}_{u|\mathbf{v}} \ell^\dagger} \right|^2. \tag{S26}$$

We proceed to bound the two terms in the equation above. The following remark will be very useful throughout this proof.

**Remark** 5. Let $(\mathcal{X}, \|\cdot\|)$ be a normed vector space. Let $a, a_h \in \mathcal{X}$, where $a_h$ depends on some discretisation parameter $h$. Suppose there are constants $\widetilde{M}, p > 0$ such that we have $\|a - a_h\| \leq \widetilde{M} h^p + \mathcal{O}(h^{p+1})$ as $h \to 0$. Then we can bound $\|a_h\|$ as follows:

$$\|a_h\| = \|a + (a_h - a)\|$$
$$\leq \|a\| + \|a_h - a\|$$
$$\leq \|a\| + \widetilde{M} h^p + \mathcal{O}(h^{p+1}) \tag{S27}$$
$$\leq \|a\| + \mathcal{O}(h^p). \tag{S28}$$

Both of the forms Eq. (S27) and Eq. (S28) will be useful in what follows. $\blacklozenge$

It will be helpful to first derive a bound on $\|\Sigma_\star - \Sigma_h\|_{L^2}$. We proceed as follows:

$$\|\Sigma_\star - \Sigma_h\|_{L^2} = \|\mathcal{L}^{-1} K \mathcal{L}^{-\dagger} - \hat{L}_h^{-1} K \hat{L}_h^{-\dagger}\|_{L^2}$$
$$\leq \|\mathcal{L}^{-1} K (\mathcal{L}^{-1} - \hat{L}_h^{-1})^\dagger\|_{L^2} + \|(\mathcal{L}^{-1} - \hat{L}_h^{-1}) K \hat{L}_h^{-\dagger}\|_{L^2}$$
$$\leq \|\mathcal{L}^{-1}\|_{L^2} \|K\|_{L^2} \|(\mathcal{L}^{-1} - \hat{L}_h^{-1})^\dagger\|_{L^2} + \|\mathcal{L}^{-1} - \hat{L}_h^{-1}\|_{L^2} \|K\|_{L^2} \|\hat{L}_h^{-\dagger}\|_{L^2}$$
$$= (\|\mathcal{L}^{-1}\|_{L^2} + \underbrace{\|\hat{L}_h^{-1}\|_{L^2}}_{(a)}) \|K\|_{L^2} \underbrace{\|\mathcal{L}^{-1} - \hat{L}_h^{-1}\|_{L^2}}_{(b)}$$
$$\leq (2\|\mathcal{L}^{-1}\|_{L^2} + Ch^2) \|K\|_{L^2} Ch^2$$
$$\leq C_1 h^2 + \mathcal{O}(h^4) \tag{S29}$$

6

where $C_1 := 2\|\mathcal{L}^{-1}\|_{L^2}\|K\|_{L^2}C > 0$ is a constant independent of $h$. Note we have utilised Remark 5 to bound $(a)$ and Lemma 1 to bound $(b)$. We have also made use of the fact that $\|L\|_{\mathcal{X} \to \mathcal{Y}} = \|L^\dagger\|_{\mathcal{Y} \to \mathcal{X}}$ for a bounded linear operator between two Hilbert spaces. With this result we can now begin to bound the two terms in Eq. (S26). Starting with the mean term we have,

$$\left|\ell m_{u|\mathbf{v}}^{(\star)} - \ell m_{u|\mathbf{v}}^{(h)}\right| = \left|\ell(m_{u|\mathbf{v}}^{(\star)} - m_{u|\mathbf{v}}^{(h)})\right|$$
$$\leq \|\ell\|_{L^2(D) \to \mathbb{R}} \|m_{u|\mathbf{v}}^{(\star)} - m_{u|\mathbf{v}}^{(h)}\|_{L^2(D)}.$$

Using the definitions of the posterior means we can bound,

$$\|m_{u|\mathbf{v}}^{(\star)} - m_{u|\mathbf{v}}^{(h)}\|_{L^2(D)} = \|m_\star - m_h + \Sigma_\star S^\dagger B_{\epsilon,\star}^{-1}(\mathbf{v} - Sm_\star) - \Sigma_h S^\dagger B_{\epsilon,h}^{-1}(\mathbf{v} - Sm_h)\|_{L^2(D)}$$
$$\leq \underbrace{\|m_\star - m_h\|_{L^2(D)}}_{(c)} + \underbrace{\|\Sigma_\star S^\dagger B_{\epsilon,\star}^{-1} - \Sigma_h S^\dagger B_{\epsilon,h}^{-1}\|_{\mathbb{R}^s \to L^2(D)}}_{(d)} \|\mathbf{v}\|$$
$$+ \underbrace{\|\Sigma_\star S^\dagger B_{\epsilon,\star}^{-1} Sm_\star - \Sigma_h S^\dagger B_{\epsilon,h}^{-1} Sm_h\|_{L^2(D)}}_{(e)} \quad \text{(S30)}$$

where we have defined $B_{\epsilon,i} := (\epsilon^2 I + S\Sigma_i S^\dagger)$, $i \in \{\star, h\}$. The term $(c)$ was bounded by $Ch^2 \|\bar{f}\|_{L^2(D)}$ in the proof of Theorem 3. We now bound terms $(d)$ and $(e)$. For $(d)$ we have:

$$\|\Sigma_\star S^\dagger B_{\epsilon,\star}^{-1} - \Sigma_h S^\dagger B_{\epsilon,h}^{-1}\|_{\mathbb{R}^s \to L^2(D)} \leq \|(\Sigma_\star - \Sigma_h)S^\dagger B_{\epsilon,\star}^{-1}\|_{\mathbb{R}^s \to L^2(D)}$$
$$+ \|\Sigma_h S^\dagger (B_{\epsilon,\star}^{-1} - B_{\epsilon,h}^{-1})\|_{\mathbb{R}^s \to L^2(D)}$$
$$\leq \|\Sigma_\star - \Sigma_h\|_{L^2} \|S\|_{C(D) \to \mathbb{R}^s} \|B_{\epsilon,\star}^{-1}\|_{\mathbb{R}^s \to \mathbb{R}^s}$$
$$\underbrace{\|\Sigma_h\|_{L^2}}_{(f)} \|S\|_{C(D) \to \mathbb{R}^s} \underbrace{\|B_{\epsilon,\star}^{-1} - B_{\epsilon,h}^{-1}\|_{\mathbb{R}^s \to \mathbb{R}^s}}_{(g)} \quad \text{(S31)}$$

Using Remark 5 we can bound $(f)$ by $\|\Sigma_\star\|_{L^2} + \mathcal{O}(h^2)$. In order to obtain a bound for $(g)$ we first compute:

$$\|B_{\epsilon,\star} - B_{\epsilon,h}\|_{\mathbb{R}^s \to \mathbb{R}^s} = \|S(\Sigma_\star - \Sigma_h)S^\dagger\|_{\mathbb{R}^s \to \mathbb{R}^s}$$
$$\leq \|S\|_{C(D) \to \mathbb{R}^s}^2 \|\Sigma_\star - \Sigma_h\|_{L^2}$$
$$\leq C_1 \|S\|_{C(D) \to \mathbb{R}^s}^2 h^2 + \mathcal{O}(h^4). \quad \text{(S32)}$$

Since $B_{\epsilon,\star}$ and $B_{\epsilon,h}$ are invertible we have for sufficiently small $h$ that $\|B_{\epsilon,\star} - B_{\epsilon,h}\|_{\mathbb{R}^s \to \mathbb{R}^s} < 1/\|B_{\epsilon,\star}^{-1}\|_{\mathbb{R}^s \to \mathbb{R}^s}$ since Eq. (S32) gives us that $\|B_{\epsilon,\star} - B_{\epsilon,h}\|_{\mathbb{R}^s \to \mathbb{R}^s} \leq \mathcal{O}(h^2)$. Thus, we can utilise [S3, Corollary 8.2] to deduce that:

$$\|B_{\epsilon,\star}^{-1} - B_{\epsilon,h}^{-1}\|_{\mathbb{R}^s \to \mathbb{R}^s} \leq \frac{\|B_{\epsilon,\star}^{-1}\|_{\mathbb{R}^s \to \mathbb{R}^s}^2 \|B_{\epsilon,\star} - B_{\epsilon,h}\|_{\mathbb{R}^s \to \mathbb{R}^s}}{1 - \|B_{\epsilon,\star}^{-1}\|_{\mathbb{R}^s \to \mathbb{R}^s} \|B_{\epsilon,\star} - B_{\epsilon,h}\|_{\mathbb{R}^s \to \mathbb{R}^s}} \leq C_2 h^2 + \mathcal{O}(h^4) \quad \text{(S33)}$$

where $C_2 := \|B_{\epsilon,\star}^{-1}\|_{\mathbb{R}^s \to \mathbb{R}^s}^2 C_1 \|S\|_{C(D) \to \mathbb{R}^s}^2 > 0$ is a constant independent of $h$.

Returning to Eq. (S31),

$$\|\Sigma_\star S^\dagger B_{\epsilon,\star}^{-1} - \Sigma_h S^\dagger B_{\epsilon,h}^{-1}\|_{\mathbb{R}^s \to L^2(D)} \leq (C_1 h^2 + \mathcal{O}(h^4))\|S\|_{C(D) \to \mathbb{R}^s}\|B_{\epsilon,\star}^{-1}\|_{\mathbb{R}^s \to \mathbb{R}^s}$$
$$+ (\|\Sigma_\star\|_{L^2} + \mathcal{O}(h^2))\|S\|_{C(D) \to \mathbb{R}^s}(C_2 h^2 + \mathcal{O}(h^4))$$
$$\leq C_3 h^2 + \mathcal{O}(h^4) \qquad (S34)$$

where $C_3 := C_1\|S\|_{C(D) \to \mathbb{R}^s}\|B_{\epsilon,\star}^{-1}\|_{\mathbb{R}^s \to \mathbb{R}^s} + C_2\|S\|_{C(D) \to \mathbb{R}^s}\|\Sigma_\star\|_{L^2} > 0$ is a constant independent of $h$. Now returning to Eq. (S30), for $(e)$ we have:

$$\|\Sigma_\star S^\dagger B_{\epsilon,\star}^{-1} S m_\star - \Sigma_h S^\dagger B_{\epsilon,h}^{-1} S m_h\|_{L^2(D)}$$
$$\leq \|\Sigma_\star S^\dagger B_{\epsilon,\star}^{-1} S m_\star - \Sigma_\star S^\dagger B_{\epsilon,\star}^{-1} S m_h\|_{L^2(D)} + \|\Sigma_\star S^\dagger B_{\epsilon,\star}^{-1} S m_h - \Sigma_h S^\dagger B_{\epsilon,h}^{-1} S m_h\|_{L^2(D)}$$
$$= \|\Sigma_\star S^\dagger B_{\epsilon,\star}^{-1} S(m_\star - m_h)\|_{L^2(D)} + \|(\Sigma_\star S^\dagger B_{\epsilon,\star}^{-1} - \Sigma_h S^\dagger B_{\epsilon,h}^{-1})S m_h\|_{L^2(D)}$$
$$\leq \|\Sigma_\star\|_{L^2}\|S\|_{C(D) \to \mathbb{R}^s}^2\|B_{\epsilon,\star}^{-1}\|_{\mathbb{R}^s \to \mathbb{R}^s}\|m_\star - m_h\|_{L^2(D)}$$
$$+ \underbrace{\|\Sigma_\star S^\dagger B_{\epsilon,\star}^{-1} - \Sigma_h S^\dagger B_{\epsilon,h}^{-1}\|_{\mathbb{R}^s \to L^2(D)}}_{\text{bound using Eq. (S34)}} \|S\|_{C(D) \to \mathbb{R}^s} \underbrace{\|m_h\|_{L^2(D)}}_{\text{Remark 5}}$$
$$\leq \|\Sigma_\star\|_{L^2}\|S\|_{C(D) \to \mathbb{R}^s}^2\|B_{\epsilon,\star}^{-1}\|_{\mathbb{R}^s \to \mathbb{R}^s} C h^2 \|\bar{f}\|_{L^2(D)}$$
$$+ (C_3 h^2 + \mathcal{O}(h^4))\|S\|_{C(D) \to \mathbb{R}^s}(C h^2\|\bar{f}\|_{L^2(D)} + \|m_\star\|_{L^2(D)})$$
$$\leq C_4 h^2 + \mathcal{O}(h^4) \qquad (S35)$$

where $C_4 := C\|\Sigma_\star\|_{L^2}\|S\|_{C(D) \to \mathbb{R}^s}^2\|B_{\epsilon,\star}^{-1}\|_{\mathbb{R}^s \to \mathbb{R}^s}\|\bar{f}\|_{L^2(D)} + C_3\|S\|_{C(D) \to \mathbb{R}^s}\|m_\star\|_{L^2(D)} > 0$ is a constant independent of $h$. Putting the three bounds together we obtain:

$$\|m_{u|\mathbf{v}}^{(\star)} - m_{u|\mathbf{v}}^{(h)}\|_{L^2(D)} \leq C\|\bar{f}\|_{L^2(D)} h^2 + C_3\|\mathbf{v}\|h^2 + C_4 h^2 + \mathcal{O}(h^4)$$
$$= C_5 h^2 + \mathcal{O}(h^4)$$

where $C_5 := C\|\bar{f}\|_{L^2(D)} + C_3\|\mathbf{v}\| + C_4 > 0$ is a constant independent of $h$. Thus, the mean term in Eq. (S26) can be bounded, for sufficiently small $h$, as:

$$\left|\ell m_{u|\mathbf{v}}^{(\star)} - \ell m_{u|\mathbf{v}}^{(h)}\right| \leq C_5\|\ell\|_{L^2(D) \to \mathbb{R}} h^2 + \mathcal{O}(h^4). \qquad (S36)$$

For the covariance term in Eq. (S26) it will be useful to make use of the following identity which holds for all $x, y > 0$:

$$|\sqrt{x} - \sqrt{y}| = \frac{|x - y|}{\sqrt{x} + \sqrt{y}}. \qquad (S37)$$

This identity follows from writing $x - y = (\sqrt{x} - \sqrt{y})(\sqrt{x} + \sqrt{y})$.

Using identity Eq. (S37) we have:

$$\left|\sqrt{\ell\Sigma_{u|\mathbf{v}}^{(\star)}\ell^\dagger} - \sqrt{\ell\Sigma_{u|\mathbf{v}}^{(h)}\ell^\dagger}\right| = \frac{\left|\ell\Sigma_{u|\mathbf{v}}^{(\star)}\ell^\dagger - \ell\Sigma_{u|\mathbf{v}}^{(h)}\ell^\dagger\right|}{\left|\sqrt{\ell\Sigma_{u|\mathbf{v}}^{(\star)}\ell^\dagger} + \sqrt{\ell\Sigma_{u|\mathbf{v}}^{(h)}\ell^\dagger}\right|}$$

$$= \frac{\left|\ell(\Sigma_{u|\mathbf{v}}^{(\star)} - \Sigma_{u|\mathbf{v}}^{(h)})\ell^\dagger\right|}{\left|\sqrt{\ell\Sigma_{u|\mathbf{v}}^{(\star)}\ell^\dagger} + \sqrt{\ell\Sigma_{u|\mathbf{v}}^{(h)}\ell^\dagger}\right|}$$

$$\leq \frac{\|\ell\|_{L^2(D)\to\mathbb{R}}^2 \|\Sigma_{u|\mathbf{v}}^{(\star)} - \Sigma_{u|\mathbf{v}}^{(h)}\|_{L^2}}{\left|\sqrt{\ell\Sigma_{u|\mathbf{v}}^{(\star)}\ell^\dagger} + \sqrt{\ell\Sigma_{u|\mathbf{v}}^{(h)}\ell^\dagger}\right|}.$$

In order to proceed with bounding the above we first bound the term involving the difference in posterior covariances as follows:

$$\|\Sigma_{u|\mathbf{v}}^{(\star)} - \Sigma_{u|\mathbf{v}}^{(h)}\|_{L^2} = \|\Sigma_\star - \Sigma_h - \Sigma_\star S^\dagger B_{\epsilon,\star}^{-1} S\Sigma_\star + \Sigma_h S^\dagger B_{\epsilon,h}^{-1} S\Sigma_h\|_{L^2}$$
$$\leq \|\Sigma_\star - \Sigma_h\|_{L^2} + \|\Sigma_\star S^\dagger B_{\epsilon,\star}^{-1} S\Sigma_\star - \Sigma_h S^\dagger B_{\epsilon,h}^{-1} S\Sigma_h\|_{L^2} \quad (S38)$$

The first term above has already been bounded so we focus on the second term and compute:

$$\|\Sigma_\star S^\dagger B_{\epsilon,\star}^{-1} S\Sigma_\star - \Sigma_h S^\dagger B_{\epsilon,h}^{-1} S\Sigma_h\|_{L^2}$$
$$\leq \|\Sigma_\star S^\dagger B_{\epsilon,\star}^{-1} S\Sigma_\star - \Sigma_h S^\dagger B_{\epsilon,h}^{-1} S\Sigma_\star\|_{L^2} + \|\Sigma_h S^\dagger B_{\epsilon,h}^{-1} S\Sigma_\star - \Sigma_h S^\dagger B_{\epsilon,h}^{-1} S\Sigma_h\|_{L^2}$$
$$\leq \underbrace{\|\Sigma_\star S^\dagger B_{\epsilon,\star}^{-1} - \Sigma_h S^\dagger B_{\epsilon,h}^{-1}\|_{\mathbb{R}^s \to L^2(D)}}_{\text{bound using }Eq.\ (S34)} \|S\|_{C(D)\to\mathbb{R}^s}\|\Sigma_\star\|_{L^2}$$
$$+ \underbrace{\|\Sigma_h\|_{L^2}}_{Remark\ 5} \|S\|_{C(D)\to\mathbb{R}^s}^2 \underbrace{\|B_{\epsilon,h}^{-1}\|_{\mathbb{R}^s\to\mathbb{R}^s}}_{Remark\ 5} \|\Sigma_\star - \Sigma_h\|_{L^2}$$
$$\leq C_3 \|S\|_{C(D)\to\mathbb{R}^s}\|\Sigma_\star\|_{L^2} h^2$$
$$+ \|S\|_{C(D)\to\mathbb{R}^s}^2 (\|\Sigma_\star\|_{L^2} + \mathcal{O}(h^2))(\|B_{\epsilon,\star}^{-1}\|_{\mathbb{R}^s\to\mathbb{R}^s} + \mathcal{O}(h^2))(C_1 h^2 + \mathcal{O}(h^2)) + \mathcal{O}(h^4)$$
$$\leq C_6 h^2 + \mathcal{O}(h^4) \quad (S39)$$

where $C_6 := C_3\|S\|_{C(D)\to\mathbb{R}^s}\|\Sigma_\star\|_{L^2} + \|S\|_{C(D)\to\mathbb{R}^s}^2 \|\Sigma_\star\|_{L^2} \|B_{\epsilon,\star}^{-1}\|_{\mathbb{R}^s\to\mathbb{R}^s} C_1 > 0$ is a constant independent of $h$. We can now bound Eq. (S38) as follows:

$$\|\Sigma_{u|\mathbf{v}}^{(\star)} - \Sigma_{u|\mathbf{v}}^{(h)}\|_{L^2} \leq C_7 h^2 + \mathcal{O}(h^4) \quad (S40)$$



\end{document}